\documentclass[a4paper,11pt]{article}
\usepackage{a4wide}
\usepackage{theorem}
\usepackage{amsmath}
\usepackage{array}
\usepackage{amssymb}
\usepackage{amsfonts}
\usepackage[french,english]{babel}
\usepackage{epsf}
\usepackage{epsfig}

\newtheorem{theo}{\indent Theorem\newline}[section]

{\theorembodyfont{\rmfamily}%
\newtheorem{rem}[theo]{\noindent Remark}}
{\theorembodyfont{\rmfamily}%
 \theoremstyle{break}%
}
\newtheorem{prop}[theo]{\indent Proposition\newline}
\newtheorem{lemma}[theo]{\indent Lemma\newline}
\newtheorem{cor}[theo]{\indent Corollary\newline}

 \def\N{{\Bbb{N}}}
\def\Z{{\Bbb{Z}}}

\def\R{{\Bbb{R}}}
\def\C{{\Bbb{C}}}

\usepackage{eufrak}

\newcommand{\mult}{\mathop{\rm mult}\nolimits}
\newcommand{\ind}{\mathop{\rm ind}\nolimits}

\newcommand{\coker}{\mathop{\rm coker}\nolimits}

\newlength{\indentation}%
\makeatletter
\newcommand\@makefntextsans[1]{%
    \parindent 0em%
    \noindent%
    \hb@xt@0em{\hss}%
    #1}
\def\footnotetextsans{%
     \@ifnextchar [\@xfootnotenextsans%
       {\@footnotetextsans}}
\def\@xfootnotenextsans[#1]{%
  \begingroup%
     \csname c@\@mpfn\endcsname #1\relax%
  \endgroup%
  \@footnotetextsans}
\long\def\@footnotetextsans#1{\insert\footins{%
    \reset@font\footnotesize%
    \interlinepenalty\interfootnotelinepenalty%
    \splittopskip\footnotesep%
    \splitmaxdepth \dp\strutbox \floatingpenalty \@MM%
    \hsize\columnwidth \@parboxrestore%
    \color@begingroup%
      \@makefntextsans{%
        \rule\z@\footnotesep\ignorespaces#1\@finalstrut\strutbox}
    \color@endgroup}}
\makeatother

\begin{document}

\cleardoublepage
\title{Invariants of real symplectic $4$-manifolds out of reducible and 
cuspidal curves}
\author{Jean-Yves Welschinger}
\maketitle

\makeatletter\renewcommand{\@makefnmark}{}\makeatother
\footnotetextsans{Keywords :  Real symplectic manifold, 
rational curve, enumerative geometry.}
\footnotetextsans{AMS Classification : 53D45, 14N35, 14N10, 14P99.}

{\bf Abstract:}

We construct invariants under deformation of real symplectic $4$-manifolds. These
invariants are
obtained by counting three different kinds of real rational $J$-holomorphic curves which
realize a given homology class and pass through a given real configuration of (the
appropriate number of) points. These curves are cuspidal curves, reducible curves
and curves with a prescribed
tangent line at some real point of the configuration. They are counted with respect 
to some sign defined by the parity of their number of isolated real double points and 
in the
case of reducible curves, with respect to some mutiplicity. In the case of the complex
projective plane equipped with its standard symplectic form and real structure, these
invariants coincide with the ones previously constructed in \cite{Wels0}, \cite{Wels1}.
This leads to a relation between the count of real rational $J$-holomorphic curves done
in \cite{Wels0}, \cite{Wels1} and the count of real rational reducible
$J$-holomorphic curves presented here.

\section*{Statement of the results}

Let $(X , \omega , c_X)$ be a {\it real symplectic $4$-manifold}, that is a triple made of
a smooth compact $4$-manifold $X$, a symplectic form $\omega$ on $X$ and an involution $c_X$ on $X$
such that $c_X^* \omega =- \omega$. The fixed 
point set of $c_X$ is called {\it the real part of $X$} and is
denoted by $\R X$. It is assumed to be non empty here so that it is a smooth lagrangian 
surface of $(X, \omega)$. We label its connected components by $(\R X)_1 , \dots ,
(\R X)_N$. Let $l \gg 1$ be an integer large enough and ${\cal J}_\omega$ be the space of 
almost complex structures of $X$ which are tamed by $\omega$ and of class $C^l$.
Let $\R {\cal J}_\omega$ be the subspace of ${\cal J}_\omega$ made of almost complex 
structures for which the involution $c_X$ is $J$-antiholomorphic. These two spaces
are separable Banach manifolds which are non empty and contractible (see \S $1.1$
of \cite{Wels1} for the real case). Let $d \in H_2 (X ; \Z)$ be
a homology class satisfying $c_1 (X) d > 1$, $c_1 (X) d \neq 4$ and $(c_X)_* d = -d$, where $c_1 (X)$ is 
the first Chern class of the symplectic $4$-manifold $(X , \omega)$.
Let $\underline{x} = (x_1 , \dots , x_{c_1 (X) d -2}) \in X^{c_1 (X) d -2}$ be a 
{\it real configuration} of $c_1 (X) d -2$ distinct points of $X$, that is an ordered 
subset of distinct points of $X$ which is globally invariant under $c_X$. For
$j \in \{ 1 , \dots , N \}$, we set $r_j = \# ( \underline{x} \cap (\R X)_j)$ and
$r = (r_1 , \dots , r_N)$, so that the $N$-tuple $r$ encodes the equivariant isotopy
class of $\underline{x}$. We will assume throughout the paper that $r \neq (0 , \dots , 0)$,
see Remark \ref{rem0}. Finally, denote by $I$ the subset of those $i \in \{ 1 , \dots,
c_1 (X) d -2 \}$ for which $x_i$ is fixed by the involution $c_X$. For each 
$i \in I$, choose a line $T_i$ in the tangent plane $T_{x_i} \R X$.
Then, for a generic choice of $J \in \R {\cal J}_\omega$, there are only finitely many
real rational $J$-holomorphic curves which realize the homology class $d$, pass
through $\underline{x}$ and are cuspidal. Moreover, these curves are all irreducible 
and have
only transversal double points as well as a unique real ordinary cusp as singularities.
Denote by ${\cal C}usp^d (J , \underline{x})$ this finite set of cuspidal curves.
Likewise, there are only finitely many
real rational $J$-holomorphic curves which realize the homology class $d$, pass
through $\underline{x}$ and are reducible. Moreover, these curves have only two 
irreducible 
components and only transversal double points as singularities. Denote by
${\cal R}ed^d (J , \underline{x})$ this finite set of reducible curves.
Note that since $I \neq \emptyset$, both irreducible 
components of such curves are real. 
Finally, there are only finitely many
real rational $J$-holomorphic curves which realize the homology class $d$, pass
through $\underline{x}$ and whose tangent line at some point $x_i$, $i \in I$, is $T_i$. 
Moreover, the
point $x_i$ having this property is then unique and these curves are all irreducible 
with only transversal double points as singularities. Denote by 
${\cal T}an^d (J , \underline{x})$ this finite set of rational curves. Note that if
$C \in {\cal C}usp^d (J , \underline{x}) \cup {\cal R}ed^d (J , \underline{x}) \cup
{\cal T}an^d (J , \underline{x})$, then all the singularities of $C$ are disjoint
from $\underline{x}$. Following \cite{Wels0}, \cite{Wels1}, we define the {\it mass} of 
$C$ and denote by $m(C)$ its number of real isolated double points. Here, a real double
point is said to be {\it isolated} when it is the local intersection of two complex 
conjugated branches, whereas it is said to be {\it non isolated} when it is
the local intersection of two real branches. If $C$ belongs to ${\cal R}ed^d (J ,
\underline{x})$ and $C_1$, $C_2$
denote its irreducible components, then we define the {\it multiplicity} of $C$,
and denote by $\mult(C)$, the number of real intersection points between $C_1$ and $C_2$,
that is the cardinality of $\R C_1 \cap \R C_2$.
We then set:
$$\Gamma^d_r (J , \underline{x}) = \sum_{C \in {\cal C}usp^d (J , \underline{x})
\cup {\cal T}an^d (J , \underline{x})}
(-1)^{m(C)} - \sum_{C \in {\cal R}ed^d (J , \underline{x})} (-1)^{m(C)}\mult(C).$$
\begin{theo}
\label{maintheo}
Let $(X, \omega , c_X)$ be a real symplectic $4$-manifold and $d \in H_2 (X ; \Z)$ be
such that $c_1 (X) d > 1$, $c_1 (X) d \neq 4$. The connected components of $\R X$ are labeled by 
$(\R X)_1 , \dots , (\R X)_N$. Let $\underline{x} \subset X$ be a real configuration of
$c_1 (X) d - 2$ distinct points, $r_j$ be the cardinality of $\underline{x} \cap
(\R X)_j$ and $r = (r_1 , \dots , r_N)$. Finally, let $J \in \R {\cal J}_\omega$
be generic enough so that the integer $\Gamma^d_r (J , \underline{x})$ is well defined.
Then, this integer $\Gamma^d_r (J , \underline{x})$ neither depends on the choice of
$J$, nor on the choice of $\underline{x}$.
\end{theo}
(The condition $c_1 (X) d \neq 4$ is to avoid appearance of multiple curves, see
Remark \ref{remmultiple}.)

From this theorem, the integer $\Gamma^d_r (J , \underline{x})$ can 
be denoted without ambiguity by $\Gamma^d_r$, and when it is not well defined,
we set $\Gamma^d_r = 0$. We then denote by
$\Gamma^d (T)$ the generating function $\sum_{r \in \N^N} \Gamma^d_r T^r \in \Z [T_1,
\dots , T_N]$, where $T^r =T_1^{r_1} \dots T_N^{r_N}$. This polynomial function is of
the same parity as $c_1 (X)d$ and each of its monomial actually only depends on one
indeterminate. It follows from Theorem \ref{maintheo} that the
function $\Gamma : d \in H_2 (X ; \Z) \mapsto \Gamma^d (T) \in \Z [T]$ only depends of
the real symplectic $4$-manifold $(X , \omega , c_X)$. Moreover, it is invariant under
deformation of this real symplectic $4$-manifold, that is if $\omega_t$ is a continuous
family of symplectic forms on $X$ for which $c_X^* \omega_t = - \omega_t$, then
this function is the same for all $(X ,\omega_t , c_X)$. As an application of
this invariant, we obtain the following lower bounds in real enumerative geometry.
\begin{cor}
\label{corlower}
Under the hypothesis of Theorem \ref{maintheo}, the integer $|\Gamma^d_r|$ provides
a lower bound for the cardinality of the weighted set
${\cal C}usp^d (J , \underline{x}) \cup {\cal R}ed^d (J , \underline{x}) \cup
{\cal T}an^d (J , \underline{x})$, independently of the choice of a generic 
$J \in \R {\cal J}_\omega$ and $\underline{x}$. $\square$
\end{cor}
The non triviality of the invariant $\Gamma^d_r$ is guaranteed by the following 
proposition, see Corollary $1.4$ of \cite{Wels3}.
\begin{prop}
\label{proprelation}
Let $(X, \omega , c_X)$ be the complex projective plane equipped with its standard
symplectic form and real structure, so that $H_2 (X ; \Z)$ is canonically isomorphic
to $\Z$. Let $r,d$ be integers satisfying $d \geq 2$ and $1 \leq r \leq 3d-2$.
Then $\Gamma^d_r = \chi^d_{r+1}$. $\square$
\end{prop}
Remember that the integer $\chi^d_{r+1}$ has been defined in \cite{Wels0}, \cite{Wels1}
by counting the number of real rational $J$-holomorphic curves of degree $d$ which
pass through $3d-1$ points with respect to the parity of their mass. The exact value
of this invariant is only known up to degree five in $\C P^2$, see \cite{IKS2}. In
particular, no recurrence formula analogous to the one obtained by Kontsevich to
compute the rational Gromov-Witten invariants of $\C P^2$ is known. The equality
given by Proposition \ref{proprelation} provides a relation between this invariant
$\chi^d_{r+1}$ and an analogous sum over all real reducible curves passing through 
$3d-2$ points -the middle term in the expression of $\Gamma^d_r$-. That is precisely
what one would need to provide such a recurrence formula. However, the reducible curves
are counted here with respect to some real multiplicity which is not under control, 
so as the two other terms in the expression of $\Gamma^d_r$.

The paper is organized as follows. The first paragraph is devoted to the construction
of the moduli space $\R {\cal M}^d_{cusp}$ of real rational cuspidal pseudo-holomorphic curves which
realize the homology class $d$. This space is equipped with a projection
$\pi_\R : \R {\cal M}^d_{cusp} \to \R {\cal J}_\omega \times \R_\tau X^{c_1 (X)d - 2}$.
The critical points of $\pi_\R$ as well as its lack of properness are discussed there.
The second paragraph is entirely  devoted to the study of one particular type of 
critical points of $\pi_\R$, namely those arising from curves having a degenerated
cuspidal point. The third paragraph is devoted to the study of the Gromov compactification
$\overline{\R {\cal M}}^d_{cusp}$ of $\R {\cal M}^d_{cusp}$. Finally, the fourth
paragraph is devoted to the proof of Theorem \ref{maintheo}.\\

{\bf Acknowledgements:}

This work was initiated during my stay at the Mathematical Sciences Research Institute
in spring $2004$. I would like to acknowledge MSRI for the excellent 
working conditions it provided to me.

\tableofcontents

\section{Moduli space of real rational cuspidal pseudo-holomorphic curves}
\label{sectmoduli} 

Let $d \in H_2 (X ; \Z)$ be such that $(c_X)_* d = -d$ and $c_1 (X) d > 1$, $c_1 (X) d \neq 4$. Let $\tau$
be an order two permutation of the set $\{ 1 , \dots , c_1 (X) d -2 \}$ and
$c_\tau : (x_1 , \dots , x_{c_1 (X) d -2}) \in X^{c_1 (X) d -2} \mapsto
(c_X (x_{\tau (1)}) , \dots , c_X (x_{\tau (c_1 (X) d -2)}) \in X^{c_1 (X) d -2}$
be the associated real structure of $X^{c_1 (X) d -2}$. The fixed point set of
$c_{\tau}$ is denoted by $\R_\tau X^{c_1 (X) d -2}$.

\subsection{Moduli space ${\cal P}^*_{cusp}$ of cuspidal pseudo-holomorphic maps}
\label{subsectmodulicusp}

Let $S$ be an oriented sphere of dimension two and ${\cal J}_S$ be the space of complex 
structures of class $C^l$ of $S$ which are compatible with its orientation. Let
$\underline{z} = (z_1 , \dots , z_{c_1 (X) d -2}) \in S^{c_1 (X) d -2}$ be an ordered 
set of $c_1 (X) d -2$ distinct points of $S$. Let $\nabla$ be a torsion free connection
on $TX$ which is invariant under $c_X$. We set
$${\cal P} = \{ (u , J_S , J , \underline{x}) \in L^{k,p} (S,X) \times {\cal J}_S
\times {\cal J}_\omega \times  X^{c_1 (X) d -2} \, | \, u_* [S] = d \, , \,
u(\underline{z}) = \underline{x} \, , \, du + J \circ du \circ J_S = 0 \},
$$
where $1 \ll k \ll l$ is large enough and $p > 2$. Let 
${\cal P}^* \subset {\cal P}$ be the
space of {\it non multiple} pseudo-holomorphic maps, that is the space of quadruples
$(u , J_S , J , \underline{x})$ for which $u$ cannot be written $u' \circ \Phi$
where $\Phi : S \to S'$ is a non trivial ramified covering and $u' : S' \to X$
a pseudo-holomorphic map. Remember that ${\cal P}^*$ is a separable Banach manifold
of class $C^{l-k}$ (see \cite{McDSal}, Proposition $3.2.1$) with tangent bundle
$$T_{(u , J_S , J , \underline{x})} {\cal P}^* = \{ (v , \dot{J}_S ,
\dot{J} , \stackrel{.}{\underline{x}}) \in T_{(u , J_S , J , \underline{x})}
(L^{k,p} (S,X) \times {\cal J}_S
\times {\cal J}_\omega \times  X^{c_1 (X) d -2}) \, | \, v(\underline{z}) = 
\stackrel{.}{\underline{x}} $$ $$\text{and } 
Dv + J \circ du \circ \dot{J}_S +
\dot{J} \circ du \circ J_S = 0 \}.$$
Here, $T_u L^{k,p} (S,X) = \{ v \in L^{k,p} (S,E_u) \}$ where $E_u = u^* TX$ and
$D : v \in L^{k,p} (S,E_u) \mapsto \nabla v + J \circ \nabla v \circ  J_S +
\nabla_v J \circ du \circ  J_S \in L^{k-1,p} (S, \Lambda^{0,1} S \otimes E_u)$ is
the associated Gromov operator (see \cite{McDSal}, Proposition $3.1.1$).
Let 
$${\cal P}^*_{cusp} = \{ ((u , J_S , J , \underline{x}), z_c) \in {\cal P}^* \times S
\, | \, d_{z_c} u = 0 \}, \text{ and}$$
$${\cal P}^*_{hocusp} = \{ ((u , J_S , J , \underline{x}), z_c) \in {\cal P}^*_{cusp}
\, | \, \nabla du|_{z_c} = 0 \}$$
be the subspace of maps having a higher order cuspidal point at $z_c$.

\begin{prop}
\label{propp*cusp}
The space ${\cal P}^*_{cusp}$ is a separable Banach manifold of class $C^{l-k}$ with
tangent bundle $T_{((u , J_S , J , \underline{x}), z_c)} {\cal P}^*_{cusp} =
\{(v,\dot{J}_S , \dot{J} , \stackrel{.}{\underline{x}} ,
\stackrel{.}{z}_c) \in T_{(u , J_S , J , \underline{x})} {\cal P}^* 
\times T_{z_c} S \, | \, \nabla v|_{z_c} + \nabla_{\stackrel{.}{z}_c} du = 0 \}.$
The space ${\cal P}^*_{hocusp}$ is a separable Banach submanifold of ${\cal P}^*_{cusp}$ 
of class $C^{l-k}$ and real codimension four.
\end{prop}

{\bf Proof:}

The proof is analogous to the one of Proposition $2.7$ of \cite{Wels1}, we just recall
a sketch of it. Denote by $F$ the vector bundle over ${\cal P}^* \times S$ whose fibre
over $((u , J_S , J , \underline{x}), z_c)$ is the vector space $T^*_{z_c} S \otimes
T_{u(z_c)} X$. In particular, the restriction of $F$ over $\{(u , J_S , J , \underline{x})
\} \times S$ is the bundle $T^* S \otimes_\C u^* TX$. From Proposition $3.2.1$ of
\cite{McDSal}, the bundle $F$ is of class $C^{l-k}$ since trivialization maps depend
$C^{l-k}$-smoothly on $u$ and $C^{l-2}$-smoothly on $z_c$, $u$ being of class $C^l$
from \cite{McDSal}, Theorem $B.4.1$. The section $d_{z_c} u$ of $F$ is of class
$C^{l-k}$ and vanishes transversely from Lemma $2.6$ of \cite{Wels1}. The first part
of the proposition follows and the second part can be proved along the same lines.
$\square$\\

Remember that if $o(z_c)$ denotes the vanishing order of $du$ at $z_c$, then the
jet of $u$ at the order $2 o(z_c) + 1$ is a well defined complex polynomial
(see \cite{Sik}, Proposition $3$). The subspace ${\cal P}^*_{hocusp}$ is precisely
made of maps $u$ for which $o(z_c) > 1$. When $o(z_c)=1$, this complex polynomial
can be written $j_2 (u) (z - z_c)^2 + j_3 (u) (z - z_c)^3$ where $j_2 (u) , j_3 (u)
\in T_{u(z_c)} X$, $z$ is a complex coordinate of $(S , J_S)$ in a neighbourhood of
$z_c$ and $j_2 (u) \neq 0$ generates the tangent line of $u$ at the cuspidal point
$u(z_c)$. The cuspidal points for which $j_3 (u)$ is colinear to $j_2 (u)$ are
said to be {\it degenerated}. They will be studied in detail in \S 
\ref{sectiondegcusp}.

\subsection{Normal sheaf}
\label{subsectionnormal}

Remember that the $\C$-linear part of the Gromov operator $D$ 
is some $\overline{\partial}$-operator denoted by $\overline{\partial}$. The latter induces a
holomorphic structure on the bundle $E_u = u^* TX$ which turns the morphism $du : TS \to
E_u$ into an injective homomorphism of analytic sheaves (see \cite{IShev}, Lemma $1.3.1$).
Likewise, the $\C$-antilinear part of $D$ is some order $0$ operator 
denoted by $R$ and defined by the formula $R_{(u,J_S , J, \underline{x})} (v) = 
N_J (v, du)$ where $N_J$ is the Nijenhuis tensor of $J$. Denote by 
${\cal N}_u$ the quotient sheaf ${\cal O}_S (E_u)/du({\cal O}_S (TS))$ so that it fits in
the following exact sequence of analytic sheaves $0 \to {\cal O}_S (TS)\to 
{\cal O}_S (E_u) \to {\cal N}_u \to 0$. As soon as $((u,J_S , J, \underline{x}), z_c)$
belongs to ${\cal P}^*_{cusp}$, this exact sequence extends to
$0 \to {\cal O}_S (TS) \otimes {\cal O}_S (z_c) \to 
{\cal O}_S (E_u) \to {\cal N}_u^{z_c} \to 0$, where ${\cal N}_u^{z_c}$ is a quotient
of the sheaf ${\cal N}_u$. We denote in this case by $\C_{z_c}$ the skyscraper
subsheaf $du({\cal O}_S (TS) \otimes {\cal O}_S (z_c)) / du({\cal O}_S (TS))$ of
${\cal N}_u$. Also, in this case, we denote by $E_u^{cusp}$ the subsheaf
$\{ v \in {\cal O}_S (E_u) \, | \, \nabla v|_{z_c} \in Im (\nabla du|_{z_c}) \}$
of $E_u$ and by ${\cal N}_u^{cusp}$ the quotient sheaf ${\cal O}_S (E_u^{cusp}) /
du ({\cal O}_S (TS))$. We hence obtain the exact sequence
$0 \to {\cal O}_S (TS) \to {\cal O}_S (E_u^{cusp}) \to {\cal N}_u^{cusp} \to 0$.
Note that from the inclusion ${\cal O}_S (E_u^{cusp}) \subset {\cal O}_S (E_u)$
follows the inclusion ${\cal N}_u^{cusp} \subset {\cal N}_u^{z_c} \oplus \C_{z_c}$.
Finally, denote by ${\cal O}_S (TS_{- \underline{z}})$ (resp. ${\cal O}_S (E_{u , 
- \underline{z}})$, ${\cal O}_S (E^{cusp}_{u , 
- \underline{z}})$, ${\cal N}_{u , 
- \underline{z}}$, ${\cal N}^{z_c}_{u , 
- \underline{z}}$, ${\cal N}^{cusp}_{u , 
- \underline{z}}$, $\C_{z_c, 
- \underline{z}}$) the subsheaf of sections of ${\cal O}_S (TS)$ 
(resp. ${\cal O}_S (E_{u})$, ${\cal O}_S (E^{cusp}_{u})$, ${\cal N}_{u}$, 
${\cal N}^{z_c}_{u}$, ${\cal N}^{cusp}_{u}$ , $\C_{z_c}$) which vanish at 
$\underline{z}$.

\begin{prop}
\label{propprepcrit}
Let $((u,J_S , J, \underline{x}) , z_c) \in {\cal P}^*_{cusp} \setminus 
{\cal P}^*_{hocusp}$. Then, we have the inclusion $\C_{z_c , 
- \underline{z}} \subset {\cal N}^{cusp}_{u , 
- \underline{z}}$ if and only if $z_c \in \underline{z}$ or $z_c$ is a degenerated
cuspidal point of $u$. In both cases, $\nabla$ induces at $z_c$ a derivation
$\nabla^{z_c}$ of sections of ${\cal N}^{z_c}_{u , 
- \underline{z}}$ such that the image of ${\cal N}^{cusp}_{u , 
- \underline{z}}$ in ${\cal N}^{z_c}_{u , 
- \underline{z}}$ under the projection ${\cal N}_{u , 
- \underline{z}} \to {\cal N}^{z_c}_{u , 
- \underline{z}}$ with kernel $\C_{z_c, 
- \underline{z}}$ is the subsheaf $\{ v \in {\cal N}^{z_c}_{u , 
- \underline{z}} \, | \, \nabla^{z_c} v = 0 \}$.
\end{prop}
Note that if $((u,J_S , J, \underline{x}) , z_c) \in {\cal P}^*_{cusp} \setminus 
{\cal P}^*_{hocusp}$ does not satisfy $\C_{z_c , 
- \underline{z}} \subset {\cal N}^{cusp}_{u , 
- \underline{z}}$, then the projection ${\cal N}_{u , 
- \underline{z}} \to {\cal N}^{z_c}_{u , 
- \underline{z}}$ with kernel $\C_{z_c, 
- \underline{z}}$ establishes an isomorphism between the sheaves ${\cal N}^{cusp}_{u , 
- \underline{z}}$ and ${\cal N}^{z_c}_{u , - \underline{z}}$.\\

{\bf Proof:}

There exist a complex coordinate $z$ of $(S , J_S)$ in a neighbourhood $U$ of
$z_c$ as well as a local chart of $X$ in a neighbourhood of
$u(z_c)$ such that the map $u$ writes $z \in U \mapsto ( (z - z_c)^2 , a (z-z_c)^3)
+ (z-z_c)^3 \epsilon_1 (z-z_c) \in \C^2$, where $a \in \C$ and $\epsilon_1 \in 
L^{k,p} (U , \C^2)$, $\epsilon_1 (z_c) = 0$ (see \cite{Sik}, Proposition $3$).
We can assume that the connection $\nabla$ is the standard connection given by 
this chart. Then, the image $Im (\nabla du|_{z_c})$ is carried by the first coordinate axis
of $\C^2$. Now, if $z_c \notin \underline{z}$ (resp. $z_c \in  \underline{z}$),
a local section of $\C_{z_c , -\underline{z}}$ writes $v_{z_c} = du (\frac{1}{z-z_c})
= (2 , 3a(z-z_c)) + (z-z_c)\epsilon_2 (z-z_c)$ (resp. $v_{z_c} = du (1) =
(2 (z-z_c), 3a(z-z_c)^2) + (z-z_c)^2 \epsilon_2 (z-z_c)$), where $\epsilon_2 \in 
L^{k-1,p} (U , \C^2)$, $\epsilon_2 (z_c) = 0$. We deduce that $\nabla v_{z_c}|_{z_c}
\in Im (\nabla du|_{z_c})$ if and only if $a=0$ or $z_c \in \underline{z}$.
The first part of the proposition is proved. In both cases, the evaluation of $\nabla
v_{z_c}$ at the point $z_c$ vanishes in ${\cal N}^{z_c}_{u , -\underline{z}}$.
Thus, for every local section $v$ of ${\cal N}^{z_c}_{u , -\underline{z}}$, the
evaluation of $\nabla v$ at the point $z_c$ does not depend on the choice of a lift
of $v$ in ${\cal O}_S (E_{u , -\underline{z}})$. We denote by $\nabla^{z_c} v \in
T_{u(z_c)} X / Im (\nabla du|_{z_c})$ this value. A section $v$ of ${\cal O}_S 
(E_{u , -\underline{z}})$ satisfies then $\nabla v |_{z_c} \in Im (\nabla du|_{z_c})$ 
if and only if the quotient section satisfies $\nabla^{z_c} v = 0$. $\square$\\

Remember that the operator $D : L^{k,p} (S , E^{cusp}_{u , -\underline{z}}) \to
L^{k-1,p} (S , \Lambda^{0,1} S \otimes E^{cusp}_{u})$ induces a quotient
operator $\overline{D} : L^{k,p} (S , {\cal N}^{cusp}_{u , -\underline{z}}) :=
L^{k,p} (S , E^{cusp}_{u , -\underline{z}}) / du(L^{k,p} (S , TS_{-\underline{z}})
\to$ $L^{k-1,p} (S , \Lambda^{0,1} S \otimes N_{u})$. Here, 
$N_u$ denotes the normal bundle of $u$ and $N_{u , -\underline{z}} = N_u \otimes
{\cal O}_S (-\underline{z})$. From the short exact sequence of complexes
$$\begin{array}{ccccccccc}
0 & \to & L^{k,p} (S , TS_{-\underline{z}}) & \stackrel{du}{\to} & 
L^{k,p} (S , E^{cusp}_{u , -\underline{z}}) & \to & 
L^{k,p} (S , {\cal N}^{cusp}_{u , -\underline{z}})
& \to & 0 \\
&&\downarrow \overline{\partial}_S && \downarrow D && \downarrow \overline{D}&&\\
0 & \to & L^{k-1,p} (S , \Lambda^{0,1} S \otimes TS) & \stackrel{du}{\to} & 
L^{k-1,p} (S , \Lambda^{0,1} S \otimes E^{cusp}_{u}) & \to & 
L^{k-1,p} (S , \Lambda^{0,1} S \otimes N_{u})
& \to & 0, \\
\end{array}$$
we deduce the long exact sequence
$0 \to H^0 (S , TS_{-\underline{z}}) \to H^0_D (S , E^{cusp}_{u , -\underline{z}}) \to 
H^0_{\overline{D}} (S , {\cal N}^{cusp}_{u , -\underline{z}}) \to
H^1 (S , TS_{-\underline{z}}) \to H^1_D (S , E^{cusp}_{u , -\underline{z}}) \to 
H^1_{\overline{D}} (S , {\cal N}^{cusp}_{u , -\underline{z}}) \to 0$,
where $H^0_D$, $H^0_{\overline{D}}$ (resp. $H^1_D$, $H^1_{\overline{D}}$) denote
the kernels (resp. cokernels) of the operators $D$, $\overline{D}$ on the
associated sheaves. In particular,
\begin{eqnarray*}
\ind_\R (\overline{D}) & = & \ind_\R (D) - \ind_\R (\overline{\partial}_S) \\
 & = & 2 ( c_1 (X)d + 1 - 2 \# \underline{z}) - 2(3 - \# \underline{z}) \\
& = & 0.
\end{eqnarray*}

\subsection{Moduli space of real rational cuspidal pseudo-holomorphic curves}

\subsubsection{Gauge action of ${\cal D}iff (S , \protect\underline{z})$ on ${\cal P}_{cusp}^*$}
\label{subsectgauge}

Denote by ${\cal D}iff (S,z)$ the group of diffeomorphisms of class $C^{l+1}$ of $S$,
which either preserve the orientation and fix $\underline{z}$, or reverse the orientation
and induce the permutation on $\underline{z}$ associated to $\tau$.
Let ${\cal D}iff^+ (S,z)$ (resp. ${\cal D}iff^- (S,z)$) be the subgroup 
of ${\cal D}iff (S,z)$ made of orientation
preserving diffeomorphisms (resp. its complement in ${\cal D}iff (S,z)$). Let
$s_*$ be the morphism ${\cal D}iff (S,z) \to
\Z / 2\Z$ having kernel ${\cal D}iff^+ (S,z)$.
The group ${\cal D}iff (S,z)$ acts on the pair $({\cal P}^*_{cusp} , {\cal P}^*_{hocusp})$ by

$$\phi . ((u,J_S , J , \underline{x}) , z_c) = \left\{ \begin{array}{rcl}
((u \circ \phi^{-1} , (\phi^{-1})^*J_S , J , \underline{x}) , \phi (z_c)) & 
\text{if} & s_*(\phi)=+1, \\
((c_X \circ u \circ \phi^{-1} , (\phi^{-1})^*J_S , \overline{c_X}^*(J) , 
c_X (\underline{x})) ,  \phi (z_c)) & 
\text{if} & s_*(\phi)=-1,
\end{array} \right. $$
where $(\phi^{-1})^* J_S = s_*(\phi) d\phi \circ J_S \circ d\phi^{-1}$ and $\overline{c_X}^*(J) =
- dc_X \circ J \circ dc_X$.
With the exception of the identity, only the order two elements of 
${\cal D}iff^- (S,z)$ may have 
non empty fixed 
point set in ${\cal P}^*_{cusp}$. In particular, two such involutions have disjoint 
fixed point sets (compare \cite{Wels1}, Lemma $1.3$). Moreover, the operators
$D$ and $\overline{D}$ are ${\cal D}iff (S,z)$ equivariants (compare \cite{Wels1},
Lemma $1.5$). Now, if $((u,J_S , J , \underline{x}) , z_c)$ is fixed by
some order two element $c_S$ of ${\cal D}iff^- (S,z)$, we denote by
$H^0_D (S ,  E^{cusp}_{u , -\underline{z}})_{\pm 1}$, 
$H^1_D (S ,  E^{cusp}_{u , -\underline{z}})_{\pm 1}$ (resp. 
$H^0_D (S , {\cal N}^{cusp}_{u , -\underline{z}})_{\pm 1}$, 
$H^1_D (S , {\cal N}^{cusp}_{u , -\underline{z}})_{\pm 1}$) the eigenspaces associated to 
the eigenvalue $\pm 1$ of the action of $c_S$ on the kernel and the cokernel of $D$
(resp. of $\overline{D}$).

\subsubsection{Moduli spaces ${\cal M}^d_{cusp}$, $\R {\cal M}^d_{cusp}$ and the
projections $\pi$, $\pi_\R$}
\label{subsectmoduli}

Denote by ${\cal M}^d_{cusp}$ (resp. ${\cal M}^d_{hocusp}$) the quotient of 
${\cal P}^*_{cusp}$ (resp. ${\cal P}^*_{hocusp}$) by the action of ${\cal D}iff^+ (S,z)$.
The projection $\pi : ((u,J_S , J , \underline{x}) , z_c) \in {\cal P}^*_{cusp}  \mapsto (J,\underline{x}) 
\in {\cal J}_\omega \times X^{c_1 (X)d - 2}$
induces on the quotient a projection ${\cal M}^d_{cusp} \to {\cal J}_\omega \times X^{c_1 (X)d - 2}$ 
still denoted
by $\pi$. 
\begin{prop}
\label{proppi}
The space  ${\cal M}^d_{cusp}$ is a separable Banach manifold
of class $C^{l-k}$, and $\pi$ is Fredholm of vanishing index. Moreover, the space
${\cal M}^d_{hocusp}$ is a separable Banach submanifold
of class $C^{l-k}$ and complex codimension two of ${\cal M}^d_{cusp}$. Finally, if
$[(u,J_S , J , \underline{x}) , z_c] \in {\cal M}^d_{cusp} \setminus
{\cal M}^d_{hocusp}$, then we have the isomorphisms
$\ker d\pi|_{((u,J_S , J , \underline{x}) , z_c)} \cong H^0_D (S , {\cal N}^{cusp}_{u , -\underline{z}})$
and $\coker d\pi|_{((u,J_S , J , \underline{x}) , z_c)} \cong 
H^1_D (S , {\cal N}^{cusp}_{u , -\underline{z}})$.
\end{prop}

{\bf Proof:}

The proof is analogous to the one of Corollary $2.2.3$ of \cite{Shev}.
The action of ${\cal D}iff^+ (S,z)$ on ${\cal P}^*_{cusp}$ and ${\cal P}^*_{hocusp}$
is smooth, fixed point free and admits a closed supplement. From Proposition \ref{propp*cusp}
thus follows that ${\cal M}^d_{cusp}$ and ${\cal M}^d_{hocusp}$ are separable Banach manifolds,
the latter being of codimension four in the former. Moreover, since
$\nabla du|_{z_c} \neq 0$, we have
\begin{eqnarray*}
\ker d\pi|_{[(u,J_S , J , \underline{x}) , z_c]} & = &
\{ (v, \dot{J}_S , 0 , 0 , \stackrel{.}{z}_c) \in T_{(u,J_S , J , \underline{x})} {\cal P}^*
\times T_{z_c} S \, | \, \nabla v|_{z_c} + \nabla_{\stackrel{.}{z}_c} du = 0 \}/T_{Id} {\cal D}iff^+ (S,z) \\
 & = & \{ (v, \dot{J}_S) \in L^{k , p} (S , E^{cusp}_{u , -\underline{z}}) \times
T_{J_S} {\cal J}_S \, | \, Dv = -J \circ du \circ \dot{J}_S \}/T_{Id} {\cal D}iff^+ (S,z) \\
 & = & \{ v \in L^{k , p} (S , E^{cusp}_{u , -\underline{z}}) \, | \, \exists \phi \in 
L^{k-1,p} (S , \Lambda^{0,1} S \otimes TS) \, , \, D v = du (\phi) \}/ du (L^{k,p} (S , TS_{-\underline{z}})) \\
 & = & H^0_D (S , {\cal N}^{cusp}_{u , -\underline{z}}),
\end{eqnarray*}
from the long exact sequence given at the end of \S \ref{subsectionnormal}.
Likewise, 
$$Im d\pi|_{[(u,J_S , J , \underline{x}) , z_c]} =  \{ (\dot{J} , \stackrel{.}{\underline{x}})
\in T_J {\cal J}_\omega \times T_{\underline{x}} X^{c_1 (X)d - 2} \, | \, \exists (v , \dot{J}_S ,
\stackrel{.}{z}_c) \in  L^{k , p} (S , E_{u , -\underline{z}}) \times T_{J_S} {\cal J}_S $$
$$\times
T_{z_c} S \, , \,
Dv + J \circ du  \circ \dot{J}_S = - \dot{J} \circ du \circ J_S ,
\, \nabla v|_{z_c} + \nabla_{\stackrel{.}{z}_c} du = 0 , \, v(\underline{z}) = \stackrel{.}{\underline{x}} \}$$
$$\text{so that }
\coker d\pi|_{[(u,J_S , J , \underline{x}) , z_c]} =  L^{k-1 , p} (S , \Lambda^{0,1} S 
\otimes E_u) \times T_{\underline{x}} X^{c_1 (X)d -2} /
Im(\widehat{D} \times ev),$$
where $\widehat{D} : (v,\dot{J}_S) \in L^{k , p} (S , E^{cusp}_u)
\times T_{J_S} {\cal J}_S \mapsto Dv + J \circ du \circ \dot{J}_S \in 
L^{k -1, p} (S , \Lambda^{0,1} S \otimes E_u)$ and $ev : v \in L^{k , p} (S , E^{cusp}_u)
\mapsto v(\underline{z}) \in T_{\underline{x}} X^{c_1 (X)d -2}$.
By definition, $\coker D = H^1_D (S , E^{cusp}_u)$. From the short exact sequence
$0 \to E^{cusp}_{u , - \underline{z}} \to E^{cusp}_u \stackrel{ev}{\to} T_{\underline{x}} 
X^{c_1 (X)d -2} \to 0$, we deduce the long exact sequence $\to H^0_D (S , E^{cusp}_u) \to
H^0 (S ,  T_{\underline{x}} 
X^{c_1 (X)d -2}) \to H^1_D (S , E^{cusp}_{u , - \underline{z}}) \to  H^1_D (S , E^{cusp}_u) \to 0$.
Hence, the cokernel of $D \times ev$ in $L^{k-1 , p} (S , \Lambda^{0,1} S \otimes E_u) \times
T_{\underline{x}} X^{c_1 (X)d -2}$ is isomorphic to $H^1_D (S , E^{cusp}_{u , - \underline{z}})$.
From the long exact sequence given at the end of \S \ref{subsectionnormal}, we deduce that
the cokernel of $\widehat{D} \times ev$ and hence the one of $d\pi|_{[(u,J_S , J , \underline{x}) , z_c]}$
is isomorphic to $H^1_D (S , {\cal N}^{cusp}_{u , -\underline{z}})$. $\square$\\

The manifolds ${\cal M}^d_{cusp}$ and ${\cal M}^d_{hocusp}$ are equipped with an action of 
the group ${\cal D}iff (S,z) /{\cal D}iff^+  (S,z) \cong
\Z / 2 \Z$. We denote by $\R {\cal M}^d_{cusp}$ and $\R {\cal M}^d_{hocusp}$ the fixed 
point sets of these actions. The projection $\pi$ is then $\Z / 2 \Z$-equivariant as soon as
${\cal J}_\omega \times X^{c_1 (X)d -2}$ is equipped with the action $\overline{c_X}^* \times
c_\tau$, where $\overline{c_X}^* : J \in  {\cal J}_\omega \mapsto -dc_X \circ J \circ dc_X \in
{\cal J}_\omega$ and $c_\tau$ has been defined at the beginning of \S \ref{sectmoduli}.
Denote by $\pi_\R$ the induced projection $\R {\cal M}^d_{cusp} \to \R {\cal J}_\omega \times
\R_\tau X^{c_1 (X)d -2}$.
\begin{prop}
\label{proppir}
The spaces $\R {\cal M}^d_{cusp}$ and $\R {\cal M}^d_{hocusp}$ are separable Banach manifolds
of class $C^{l-k}$, the latter being of codimension two in the former. Moreover, $\pi_\R$ 
is Fredholm of vanishing index. Finally,  if
$[(u,J_S , J , \underline{x}) , z_c] \in \R {\cal M}^d_{cusp} \setminus
\R {\cal M}^d_{hocusp}$, then we have the isomorphisms
$\ker d\pi_\R |_{((u,J_S , J , \underline{x}) , z_c)} \cong H^0_D (S , 
{\cal N}^{cusp}_{u , -\underline{z}})_{+1}$
and $\coker d\pi_\R |_{((u,J_S , J , \underline{x}) , z_c)} \cong 
H^1_D (S , {\cal N}^{cusp}_{u , -\underline{z}})_{+1}$. $\square$
\end{prop}

\subsection{Critical points of $\pi_\R$}
\label{subsectcritical}

\begin{lemma}
\label{lemmacrit}
The point $[(u,J_S , J , \underline{x}) , z_c] \in \R {\cal M}^d_{cusp} \setminus
\R {\cal M}^d_{hocusp}$ is critical for $\pi_\R$ if and only if one of the
following:

1) The differential $du$ vanishes outside $z_c$.

2) The cuspidal point $z_c$ is degenerated.

3) One has $z_c \in \underline{z}$.
\end{lemma}

{\bf Proof:}

From Proposition \ref{proppir}, $[(u,J_S , J , \underline{x}) , z_c]$ is a critical point of
$\pi_\R$ if and only if $H^0_D (S , 
{\cal N}^{cusp}_{u , -\underline{z}})_{+1} \cong H^1_D (S , 
{\cal N}^{cusp}_{u , -\underline{z}})_{+1} \neq \{ 0 \}$. We deduce from Proposition 
\ref{propprepcrit} that 2) and 3) are indeed critical points of $\pi_\R$, since in these
cases $\C_{z_c , -\underline{z}} \subset {\cal N}^{cusp}_{u , -\underline{z}}$ and $ \{ 0 \}  \neq
H^0 ( S , \C_{z_c , -\underline{z}})_{+1} \subset
H^0_D (S , {\cal N}^{cusp}_{u , -\underline{z}})_{+1}$. In the same way, if $z_c$ is a real
ordinary cusp of $u$ distinct from $\underline{z}$, then from Proposition \ref{propprepcrit},
$ {\cal N}^{cusp}_{u , -\underline{z}} = {\cal N}^{z_c}_{u , -\underline{z}}$. If $du$
does not vanish outside $z_c$, then ${\cal N}^{z_c}_{u , -\underline{z}} \cong
{\cal O}_S (N_{u , -\underline{z}})$ so that $H^1_D (S , 
{\cal N}^{cusp}_{u , -\underline{z}})_{+1} = \{ 0 \}$. Otherwise, the sheaf
${\cal N}^{cusp}_{u , -\underline{z}}$ carries some skyscraper part and
$H^0_D (S , 
{\cal N}^{cusp}_{u , -\underline{z}})_{+1} \neq \{ 0 \}$, hence the result. $\square$\\

A detailed study of critical points of type 2) will be carried out in \S \ref{sectiondegcusp}.
In particular, we will prove in Lemma \ref{lemmajet} that if $z$ is a local coordinate
in a neighbourhood of $z_c$ which is {\it adapted} to $u$, that is for which the order
three jet of $u$ reads $j_2 (z_c) (z - z_c)^2 + o(|z - z_c|^3)$, then the homogeneous part
of order $5$ of its jet is a complex monomial denoted by $j_5 (z_c) (z - z_c)^5$. Also,
note that if $u$ has two distinct ordinary cusps, then $H^1_D (S , 
{\cal N}^{cusp}_{u , -\underline{z}})_{+1} \cong H^1_D (S , 
{\cal N}^{z_c}_{u , -\underline{z}})_{+1}$ is of dimension one and from Riemann-Roch duality, it is
isomorphic to $H^0_{D^*} (S , K_S \otimes
{\cal N}^{z_c}_{u , -\underline{z}})_{-1}$ (see \cite{Wels1}, Lemma $1.7$).
Finally, remember that a
{\it stratum of codimension $m \in \N$} of a separable Banach manifold $N$ is by 
definition the image of a separable  Banach manifold $M$ under a smooth Fredholm
map $f$ of Fredholm index $-m$ such that all the limits of sequences
$\phi (x_n)$, where $x_n$ is a diverging sequence in $M$, belong to a countable
union of strata of higher codimensions. 

\begin{lemma}
\label{lemmagencrit}
1) The set of points $[(u,J_S , J , \underline{x}) , z_c] \in \R {\cal M}^d_{cusp}$ for which
$u$ has two distinct real ordinary cusps and ordinary double points as singularities,
all of which being outside $\underline{x}$ and for which any generator
$\psi$ of $H^0_{D^*} (S , K_S \otimes
{\cal N}^{z_c}_{u , -\underline{z}})_{-1} = H^1_D (S , 
{\cal N}^{cusp}_{u , -\underline{z}})_{+1}^*$ does not vanish at the cusps, is a stratum 
of codimension one of $\R {\cal M}^d_{cusp}$.

2)  The set of points $[(u,J_S , J , \underline{x}) , z_c] \in \R {\cal M}^d_{cusp}$ for which
$u$ has a degenerated cusp and transversal double points as singularities, all of which
being outside $\underline{x}$, and for which $j_5 (z_c)$ is not colinear to $j_2 (z_c)$
in an adapted coordinate, is a stratum 
of codimension one of $\R {\cal M}^d_{cusp}$.

3) The set of points $[(u,J_S , J , \underline{x}) , z_c] \in \R {\cal M}^d_{cusp}$ for which
$u$ has a real ordinary cusp at $z_c \in \underline{z}$ and only ordinary double points as other 
singularities, all of which being outside $\underline{x}$, is a stratum 
of codimension one of $\R {\cal M}^d_{cusp}$.

The union of all critical points of $\pi_\R$ not listed above belong to a countable union
of strata of codimension at least two in $\R {\cal M}^d_{cusp}$.
\end{lemma}

{\bf Proof:}

The case 2 follows from Proposition \ref{propdegcusp}. All the other cases can be proved
in the same way as Propositions $2.7$ and $2.8$ of \cite{Wels1}. For the sake of concision,
these proofs are not reproduced here. $\square$

\begin{rem}
There are actually only finitely many strata occuring in Lemma \ref{lemmagencrit}, since a
pseudo-holomorphic curve which realize the given homology class $d$ may have only
finitely many different types of singularities. Moreover, as soon as $k$ and $l$ are large
enough, all these strata are images of Banach manifolds of class $C^2$ at least. We will
use only the following fact: a generic path $\gamma : [0,1] \to \R {\cal J}_\omega
\times \R_\tau X^{c_1 (X)d -2}$ avoids the image under $\pi_\R$ of every stratum of
codimension at least two.
\end{rem}

The critical points of $\pi_\R$ listed in Lemma \ref{lemmagencrit} are
said to be {\it generic}.

\begin{theo}
\label{theocrit}
The generic critical points of $\pi_\R$ are non degenerated.
\end{theo}

{\bf Proof:}

In the case of critical points of type 1 given by Lemmas \ref{lemmacrit} and
\ref{lemmagencrit}, the proof is readily the same as the one of Lemma $2.13$ of
\cite{Wels1}. It is not reproduced here. Let $[(u,J_S , J , \underline{x}) , z_c] 
\in \R {\cal M}^d_{cusp}$ be a generic critical point of type 2 or 3 given by Lemma
\ref{lemmagencrit}. We have to prove that the quadratic form $\nabla 
d\pi_\R |_{[(u,J_S , J , \underline{x}) , z_c]} : \ker d\pi_\R |_{[(u,J_S , J , \underline{x}) , z_c]}
\times \ker d\pi_\R |_{[(u,J_S , J , \underline{x}) , z_c]} \to \coker 
 d\pi_\R |_{[(u,J_S , J , \underline{x}) , z_c]}$ is non degenerated.
Write $\pi_\R = (\pi_\R^1 , \pi_\R^2)$, $\R {\cal M}^d_{cusp} (\underline{x}) =
(\pi_\R^2)^{-1} (\underline{x})$ and $\pi_\R^J : [(u,J_S , J , \underline{x}) , z_c] \in
\R {\cal M}^d_{cusp} (\underline{x}) \mapsto J \in \R {\cal J}_\omega$ the restriction
of $\pi_\R^1$ to $\R {\cal M}^d_{cusp} (\underline{x})$. The quadratic forms
$\nabla 
d\pi_\R |_{[(u,J_S , J , \underline{x}) , z_c]}$ and $\nabla 
d\pi_\R^J |_{[(u,J_S , J) , z_c]}$ are of the same nature. Moreover, the kernel and cokernel
of the map $d\pi_\R^J$ are the same as the ones of the operator $-\widehat{D}_\R :
(v , \dot{J}_S , \dot{J} , \stackrel{.}{z}_c) \in T_{[(u,J_S , J) , z_c]}
\R {\cal M}^d_{cusp} (\underline{x})  \mapsto \dot{J} \circ du \circ J_S \in 
L^{k-1 , p} (S , \Lambda^{0,1} S \otimes N_u)$. From the relation $Dv + J \circ du \circ 
\dot{J}_S + \dot{J} \circ du \circ J_S = 0$, we deduce that
$\widehat{D}_\R (v , \dot{J}_S , \dot{J} , \stackrel{.}{z}_c) = Dv +
J \circ du \circ \dot{J}_S$. We then have to prove that 
$\nabla \widehat{D}_\R|_{[(u,J_S , J) , z_c]} : 
H^0_D (S , {\cal N}_{u, -\underline{z}}^{cusp})_{+1}^2 \to
H^1_D (S , {\cal N}_{u, -\underline{z}}^{cusp})_{+1}$ is non degenerated. Let $(v_1 , 
\dot{J}_S^1 , 0 , \stackrel{.}{z}_c^1)$ be a generator of
$H^0_D (S , {\cal N}_{u, -\underline{z}}^{cusp})_{+1}$,
then from Proposition \ref{propprepcrit}, $v_1 = du(\tilde{v}_1)$ where $\tilde{v}_1 \in
L^{k,p} (S , TS_{-\underline{z}} \otimes {\cal O}_S (z_c))_{+1}$, that is $\tilde{v}_1$ is
a meromorphic vector field on $S$ either having a simple pole at $z_c$ if $z_c \notin \underline{z}$
or which does not vanish at $z_c$ otherwise. In the same way, let $(v_2 , 
\dot{J}_S^2 , 0 , \stackrel{.}{z}_c^2) \in T_{[u, J_S , J , z_c]} 
\R {\cal M}^d_{cusp} (\underline{x})$, then
$$(\nabla_{(v_2 , 
\dot{J}_S^2 , 0 , \stackrel{.}{z}_c^2)} \widehat{D}_\R )(v_1 , 
\dot{J}_S^1 , 0 , \stackrel{.}{z}_c^1) = \big( \nabla_{(v_2 , 
\dot{J}_S^2 , 0 , \stackrel{.}{z}_c^2)} D_\R \big) (v_1) + (\nabla_{v_2} du) \circ J_S \circ
\dot{J}_S^1 \mod(Im(du)).$$
Moreover, after differentiation of the relation $D \circ du = du \circ \overline{\partial}_S$,
we deduce
$$\big( \nabla_{(v_2 , 
\dot{J}_S^2 , 0 , \stackrel{.}{z}_c^2)} D_\R \big) (du(\tilde{v}_1)) + D_\R \circ
(\nabla_{v_2} du) (\tilde{v}_1) = \nabla_{v_2} du \circ \overline{\partial}_S (\tilde{v}_1)
 \mod(Im(du)).$$
Since the relation $0 = Dv_1 + J \circ du \circ \dot{J}_S^1 = du (\overline{\partial}_S (\tilde{v}_1) + J_S 
\dot{J}_S^1)$ forces
$\overline{\partial}_S (\tilde{v}_1) + J_S \dot{J}_S^1 =0$, we deduce
$$(\nabla_{(v_2 , 
\dot{J}_S^2 , 0 , \stackrel{.}{z}_c^2)} \widehat{D}_\R )(v_1 , 
\dot{J}_S^1 , 0 , \stackrel{.}{z}_c^1) = - D_\R (\nabla_{v_2} du) (\tilde{v}_1)  
\mod(Im(du)).$$
We thus have to prove that the projection $D_\R (\nabla_{\tilde{v}_1} du (\tilde{v}_1)) \in
H^1_D (S ; {\cal N}_{u, -\underline{z}}^{cusp})_{+1} \cong \R$ does not vanish as
soon as $\tilde{v}_1 \neq 0$. Now the operator $D_\R : L^{k,p} (S, {\cal N}_{u, -\underline{z}})_{+1}
\to  L^{k-1,p} (S, \Lambda^{0,1} \otimes N_{u})_{+1}$ is an isomorphism,
and $L^{k,p} (S, {\cal N}^{cusp}_{u, -\underline{z}})_{+1} = \{ v \in 
L^{k,p} (S, {\cal N}_{u, -\underline{z}})_{+1} \, | \, \nabla^{z_c} v =0 \}$ (see Proposition
\ref{propprepcrit}) is a closed hyperplane of $L^{k,p} (S, {\cal N}_{u, -\underline{z}})_{+1}$.
We hence have to prove that $v = \nabla_{\tilde{v}_1} du (\tilde{v}_1)$ does not satisfy
$\nabla^{z_c} v =0$.

The map $u$ can be written in an adapted local chart at $z_c$ as
$u(z) = j_2 (z_c) (z-z_c)^2 + \dots + j_5 (z_c) (z-z_c)^5 + o(|z-z_c|^5)$ where
$(j_2 (z_c) , j_5 (z_c))$ does form a basis of $\R^2$ if $z_c \notin \underline{z}$ or
$u(z) = j_2 (z_c) (z-z_c)^2 + j_3 (z_c) (z-z_c)^3 + o(|z-z_c|^3)$ where
$(j_2 (z_c) , j_3 (z_c))$ does form a basis of $\R^2$ if $z_c \in \underline{z}$.
We can assume that the connection $\nabla$ is the standard one given by this chart. Then in the first case,
$\tilde{v}_1 = \frac{1}{z-z_c}$ and $v = (\nabla_{\tilde{v}_1} du) (\tilde{v}_1) =
\frac{2j_2 (z_c)}{(z-z_c)^2} + 20 j_5 (z_c) (z-z_c) + o(|z-z_c|) 
= 15j_5 (z_c) (z-z_c) + o(|z-z_c|) \mod (Im (du))$. And in the second case,
$\tilde{v}_1 = 1$ and $v = (\nabla_{\tilde{v}_1} du) (\tilde{v}_1) =
2j_2 (z_c) + 6 j_3 (z_c) (z-z_c) + o(|z-z_c|)$. In the first case, $\nabla v|_{z_c} = 15 j_5 (z_c) dz$
while in the second case $\nabla v|_{z_c} = 6 j_3 (z_c) dz$. In both cases, the projection
$\nabla^{z_c} v$ of $\nabla v|_{z_c}$ in the normal bundle does not satisfy
$\nabla^{z_c} v =0$, hence the result. $\square$

\subsection{Gromov compactification $\R \overline{\cal M}^d_{cusp}$ of $\R {\cal M}^d_{cusp}$}
\label{subsectcompact}

The projection $\pi_\R : \R {\cal M}^d_{cusp} \to \R {\cal J}_\omega \times 
\R_\tau X^{c_1 (X)d - 2}$ is not proper in general. The reason for this is that there might exist
some sequence $([u_n , J_S^n , J^n , \underline{x}^n , z_c^n])_{n \in \N}$ of
$\R {\cal M}^d_{cusp}$ such that $(J^n , \underline{x}^n)$ converges to $(J^\infty , 
\underline{x}^\infty) \in \R {\cal J}_\omega \times 
\R_\tau X^{c_1 (X)d - 2}$, but the image $u_n (S)$ converges to some reducible 
$J^\infty$-holomorphic curve. From Gromov compactness theorem (see \cite{McDSal}, Theorem
5.5.5), this is the only obstruction to the properness of $\pi_\R$. More precisely, this 
theorem describes how the sequence of maps $(u_n)_{n \in \N}$ does converge. There exist
smooth disjoint loops $\alpha_1 , \dots , \alpha_k$ in $S$ and a map $u_\infty : S \to X$
which contracts the loops $\alpha_1 , \dots , \alpha_k$ and whose image is the reducible
curve in the limit. Moreover, after may be changing the parameterization of $u_n$, this
sequence converges to $u_\infty$ in $C^0$ norm on the whole $S$ and in norm $L^{k,p}$
on every compact subset of $S \setminus (\cup_{i=1}^k \alpha_i)$. In particular, we have the
following alternative.  Either the limit $z_c^\infty$ of $(z_c^n)_{n \in \N}$ does not
belong to $\cup_{i=1}^k \alpha_i$, and then the curve in the limit has a cuspidal point at
$u_\infty (z_c^\infty)$. Or the limit $z_c^\infty$ of $(z_c^n)_{n \in \N}$ does 
belong to $\cup_{i=1}^k \alpha_i$, say $\alpha_1$, and then the two irreducible components
adjacent to $\alpha_1$ intersect each other with multiplicity at least two  at
$u_\infty (z_c^\infty)$, from adjunction formula for example.
The end of this paragraph is devoted to the proof that over a generic path $\gamma : t \in 
[0,1] \mapsto (J^t , \underline{x}^t) \in \R {\cal J}_\omega \times \R_\tau X^{c_1 (X)d - 2}$, 
the only reducible curves
which satisfy one of these two conditions have two irreducible components, 
both real, and only
transversal double points as singularities with the exception of a unique real ordinary cusp
or a unique real ordinary tacnode at some intersection point between the two irreducible components. 
Moreover, all these singularities are outside $\underline{x}^t$.

Let $m_1 \in \N$, $d_1 \in H_2 (X ; \Z)$ and $\underline{z}^1 = (z_1^1 , \dots , z_{m_1}^1) \in
S^{m_1}$ be an $m_1$-tuple of distinct points of $S$. Denote by
$$\R {\cal M}^{(d_1 , m_1)} = \{ (u^1 , J_S^1 , J^1 , \underline{x}^1) \in L^{k,p} (S , X)
\times {\cal J}_S \times {\cal J}_\omega \times X^{m_1} \, | \, du^1 + J^1 \circ du^1 
\circ J_S^1 =0 $$ 
$$\text{ and } u(\underline{z}^1) = \underline{x}^1 \}/{\cal D}iff^+ (S,\underline{z}^1).$$
Let $m_2 \in \N$, $d_2 \in H_2 (X ; \Z)$ and $\underline{z}^2 = (z_1^2 , \dots , z_{m_2}^2) \in
S^{m_2}$ be an $m_2$-tuple of distinct points of $S$. We denote by
$$\R {\cal M}^{(d_1 , m_1), (d_2 , m_2)} = \big( \R {\cal M}^{(d_1 , m_1)} \times_{{\cal J}_\omega}
\R {\cal M}^{(d_2 , m_2)} \big) \setminus Diag,$$
where $Diag = \{ ((u^1 , J_S^1 , J , \underline{x}^1),(u^2 , J_S^2 , J , \underline{x}^2)) \in
\R {\cal M}^{(d_1 , m_1)} \times_{{\cal J}_\omega}
\R {\cal M}^{(d_2 , m_2)} \, | \, u_1 (S) = u_2 (S) \}.$
We recall the following proposition (see \cite{Wels1}, Proposition $2.9$ and Corollary $2.10$).

\begin{prop}
\label{propredrap}
The space $\R {\cal M}^{(d_1 , m_1), (d_2 , m_2)}$ is a separable Banach manifold of class $C^{l-k}$. Moreover,
the projection $\pi^{(d_1 , m_1), (d_2 , m_2)}_\R : \R {\cal M}^{(d_1 , m_1), (d_2 , m_2)} \to
\R {\cal J}_\omega \times \R_{\tau_1} X^{m_1} \times \R_{\tau_2} X^{m_2}$ is Fredholm of index
$c_1 (X) d - 2 -m_1 - m_2$ where $d = d_1 + d_2$. Finally, the cokernel of
$\pi^{(d_1 , m_1), (d_2 , m_2)}_\R$ at $((u^1 , J_S^1 , J , \underline{x}^1),(u^2 , J_S^2 , J , 
\underline{x}^2))$ is isomorphic to $H^1_D (S , N_{u_1 , 
-\underline{z}^1})_{+1}
\oplus H^1_D (S , N_{u_2 , -\underline{z}^2})_{+1}$. $\square$
\end{prop} 

\begin{rem}
\label{remmultiple}
Over a generic path $\gamma : t \in 
[0,1] \mapsto (J^t , \underline{x}^t) \in \R {\cal J}_\omega \times \R_\tau X^{c_1 (X)d - 2}$,
no $J^t$-holomorphic curve which realize $d$ and pass through $\underline{x}^t$ can be multiple
or come from the diagonal $Diag$. Indeed, the condition for a $J$-holomorphic curve which
realize the homology class $\frac{1}{2}d$ to pass through $c_1 (X)d - 2$ distinct points is of
codimension $\frac{1}{2}c_1 (X)d - 1$, that is of codimension greater than one as soon as
$c_1 (X)d >4$. Moreover, a generic immersed rational $J$-holomorphic curve which realize the
homology class $\frac{1}{2}d$ has $\frac{1}{2} \big( (\frac{d}{2})^2 - \frac{1}{2}c_1(X) d +2 \big)$
transversal double points. Each of these double points is responsible for four double points of the
doubled curve. The number of double points of the
doubled curve would then be at least $2\big( (\frac{d}{2})^2 - \frac{1}{2}c_1(X) d +2 \big) - 1$,
which is impossible as soon as $c_1 (X)d <4$.
\end{rem}

\begin{prop}
\label{propredgen}
The subspace of $\R {\cal M}^{(d_1 , m_1), (d_2 , m_2)}$ made of couples  
$((u^1 , J_S^1 , J , \underline{x}^1),(u^2 , J_S^2 , J , \underline{x}^2))$ for which
$u^1$ or $u^2$ has a unique cuspidal point which is ordinary, or for which $u^1 (S)$ and
$u^2 (S)$ have a unique point of contact which is of order two, all the singularities of
$u^1 (S) \cup u^2 (S)$ being outside $\underline{x}^1 \cup \underline{x}^2$ is a stratum
of codimension one. The subspace of curves having degenerated cuspidal points or higher
order cuspidal points, or points of contact of higher order is a stratum of codimension 
at least two. $\square$
\end{prop}
We denote by $\R {\cal M}^{(d_1 , m_1), (d_2 , m_2)}_{cusp}$ (resp. $\R {\cal M}^{(d_1 , m_1), 
(d_2 , m_2)}_{tac}$) the codimension one stratum of $\R {\cal M}^{(d_1 , m_1), (d_2 , m_2)}$
given by Proposition \ref{propredgen} made of curves having a real ordinary cusp (resp.
an ordinary point of contact between $u^1 (S)$ and
$u^2 (S)$).

\begin{cor}
\label{corredgen}
Let $\gamma : t \in 
[0,1] \mapsto (J^t , \underline{x}^t) \in \R {\cal J}_\omega \times \R_\tau X^{c_1 (X)d - 2}$
be a generic path. Assume that a sequence of elements of $\R {\cal M}^d_{cusp}$ over $\gamma$
converges
to some reducible curve. Then this reducible $J^t$-holomorphic curve is given by an element
$([u^1 , J_S^1 , J , \underline{x}^1],[u^2 , J_S^2 , J , \underline{x}^2]) \in
\R {\cal M}^{(d_1 , m_1), (d_2 , m_2)}_{cusp} \cup \R {\cal M}^{(d_1 , m_1), 
(d_2 , m_2)}_{tac}$ such that $d_1 + d_2 = d$, $m_1 + m_2 = c_1 (X)d - 2$ and
$t \in ]0,1[$. Moreover, either $m_1 = c_1 (X)d_1 - 1$, or $m_1 = c_1 (X)d_1 - 2$ and
then the cuspidal point, if it exists, belongs to $u^1 (S)$.
\end{cor}

{\bf Proof:}

From Proposition \ref{propredgen} follows that the curve in the limit must belong to
$\R {\cal M}^{(d_1 , m_1), (d_2 , m_2)}_{cusp} \cup \R {\cal M}^{(d_1 , m_1), 
(d_2 , m_2)}_{tac}$ as soon as $\gamma$ is generic enough.
Now, from Proposition \ref{propredrap}, the cokernel of 
$d\pi^{(d_1 , m_1), (d_2 , m_2)}_\R |_{([u^1 , J_S^1 , J , \underline{x}^1],
[u^2 , J_S^2 , J , \underline{x}^2])}$ is isomorphic to 
$H^1_D (S , N_{u_1 , -\underline{z}^1})_{+1}
\oplus H^1_D (S , N_{u_2 , -\underline{z}^2})_{+1}$. Since this cokernel 
is of dimension less than two, we have $m_i \leq c_1 (X)d_i$ for $1 \leq i \leq 2$. Moreover,
in case $m_2 = c_1 (X)d_2$, we have 
$\dim H^1_D (S , N_{u_2 , -\underline{z}^2})_{+1} \geq 1$ with equality
if and only if the map $u^2$ is an immersion (see \cite{HLS}). Finally, the relation 
$m_1 + m_2 = c_1 (X) d -2$
forces each irreducible component in the limit to be simply covered, unless $c_1 (X) d_i \leq 1$
for some $i \in \{ 1, 2 \}$. Now, when $c_1 (X) d_i \leq 1$
for some $i \in \{ 1, 2 \}$, it suffices to count the number of double points of these rational curves
as in Remark \ref{remmultiple} to see that these irreducible components cannot be multiply covered. $\square$

\begin{rem} 
\label{rem0bis}
If we would not have excluded the case $r = (0, \dots , 0)$, then a sequence of real rational cuspidal
$J$-holomorphic curves could converge to a reducible curve having two irreducible components which
are complex conjugated and transversal to each other except at one point which is of order two, that
is an ordinary tacnode.
\end{rem}

\section{Study of degenerated order two cuspidal points}
\label{sectiondegcusp}

\subsection{Local study of degenerated order two cuspidal points}
\label{subseclocalstudy}

Let $B^4$ be the unit ball of $\C^2$ and $c_X$ be the restriction of the complex conjugation
to $B^4$. Denote by $\R B^4$ the fixed point set of $c_X$, it is the unit ball of $\R^2 \subset
\C^2$. Denote by $\R {\cal J}_{\omega_{st}}$ the space of almost complex structures $J$ of
$B^4$ which are tamed by the standard symplectic form $\omega_{st}$ and for which
$c_X$ is $J$-antiholomorphic. Let $\overline{B}^2$ be the closed unit ball of $\C$ and
$conj$ the restriction of the complex conjugation to $\overline{B}^2$. Its fixed point set
is $]-1,1[ \subset \overline{B}^2$. Finally, denote by $J_{st}$ the restriction of the complex
structure of $\C^2$ to $B^4$, so that $J_{st} \in \R {\cal J}_{\omega_{st}}$. Let $\eta > 0$
and
$$\R {\cal P}'_{cusp} (\eta) = \{ (u, J, z_c) \in L^{k,p} (\overline{B}^2 , B^4) \times
\R {\cal J}_{\omega_{st}} \times ]-1 , 1[ \, | \,  ||J - J_{st}||_{C^1} < \eta \, , \,
du + J \circ du \circ i = 0 \, , $$ $$ c_X \circ u = u \circ conj \, , \, d_{z_c} u = 0
\text{ but } \nabla du|_{z_c} \neq 0 \text{ and } u(\overline{B}^2) \text{ has smooth
 boundary} \}.$$
In particular, $u$ is not multiple. Note that $\R {\cal P}'_{cusp} (\eta)$ is not connected.
Indeed, two disks which do not have the same number of double points cannot be in the same
connected component. We are in fact interested here in a connected component for which
general elements are disks with one ordinary cusp at $z_c$ and one transversal double point.

\begin{lemma}
\label{lemmarp'}
As soon as $\eta$ is small enough, $\R {\cal P}'_{cusp} (\eta)$ is a separable Banach manifold
of class $C^{l-k}$ whose tangent space at $(u , J , z_c)$ is 
$$T_{(u , J , z_c)} \R {\cal P}'_{cusp} (\eta) = \{ (v , \dot{J} , \stackrel{.}{z}_c) \in
L^{k,p} (\overline{B}^2 , \C^2) \times T_J \R {\cal J}_{\omega_{st}} \times \R \, | \, 
Dv+\dot{J} \circ du \circ i = 0 \, , $$ $$ v = dc_X \circ v \circ conj \, , \,
\nabla v|_{z_c} + \nabla_{\stackrel{.}{z}_c} du = 0 \}. \quad \square$$
\end{lemma}

Let $(u, J, z_c) \in \R {\cal P}'_{cusp} (\eta)$, then the order three jet of $u$ at the point
$z_c$ is a complex polynomial which can be written $u(z_c) + j_2 (z_c) (z-z_c)^2 +
j_3 (z_c) (z-z_c)^3$ with $0 \neq j_2 (z_c) = \frac{1}{2} \frac{\partial^2 u}{\partial z^2} \in
\R^2$ and $j_3 (z_c) = \frac{1}{6} \frac{\partial^3 u}{\partial z^3} \in \R^2$ (see \cite{Sik},
Proposition $3$). The cuspidal point $z_c$ is {\it degenerated} when $j_3 (z_c)$ is colinear
to $j_2 (z_c)$. Let $F_{deg} : (u, J, z_c) \in \R {\cal P}'_{cusp} (\eta) \mapsto
\det(j_2 (z_c) , j_3 (z_c)) \in \R$, and $\R {\cal P}'_{dcusp} (\eta) = F_{deg}^{-1} (0)$ be
the locus of curves $(u, J, z_c)$ having a degenerated cuspidal point at $z_c$.

\begin{lemma}
\label{lemmarp'd}
As soon as $\eta$ is small enough, $\R {\cal P}'_{dcusp} (\eta)$ is a separable Banach submanifold
of class $C^{l-k}$ of $\R {\cal P}'_{cusp} (\eta)$ whose tangent space at $(u , J , z_c)$ is 
$$T_{(u , J , z_c)} \R {\cal P}'_{dcusp} (\eta) = \{ (v , \dot{J} , \stackrel{.}{z}_c) \in
L^{k,p} (\overline{B}^2 , \C^2) \times T_J \R {\cal J}_{\omega_{st}} \times \R \, | \, 
\det(d_{(v , \dot{J} , \stackrel{.}{z}_c)} j_2 (z_c) , j_3 (z_c)) $$ $$ +
\det(j_2 (z_c) ,d_{(v , \dot{J} , \stackrel{.}{z}_c)} j_3 (z_c)) = 0 \}.$$
\end{lemma}

{\bf Proof:}

The function $F_{deg}$ is of class $C^{l-k}$ as soon as $k \geq 3$ since $u$ is of class $C^l$.
It suffices to prove that $d_{(u, J, z_c)} F_{deg}$ is surjective at each point $(u, J, z_c)$
of $\R {\cal P}'_{dcusp} (\eta)$. Let $(u_0, J_0, z_c)$ be such a point. From Lemma $2.5$ of
\cite{Wels1}, there exists a smooth family of $J_0$-holomorphic maps 
$(\tilde{u}_\lambda)_{\lambda \in ]-\epsilon , \epsilon[}$  defined in a neighbourhood $U$ of
$z_c$ by $\tilde{u}_\lambda (z) = u_0 (z) + (z - z_c)^3(\lambda w + w_\lambda (z))$, where
$w$ can be any vector in $\R^2$ and $w_\lambda \in L^{k,p} (\overline{B}^2 , \C^2)$ is real
and satisfies $w_0 = 0$, $\frac{d}{d \lambda} w_\lambda|_{\lambda = 0} =0$. This family can be
extended to a smooth family 
$(u_\lambda, J_\lambda, z_c)_{\lambda \in ]-\epsilon , \epsilon[} \in 
\R {\cal P}'_{cusp} (\eta)$ such that $J_\lambda = J_0$ if $\lambda = 0$ and $J_\lambda$ differs
from $J_0$ only in a neighbourhood of $u_0 (\partial U)$. Indeed, it suffices to glue the
map $\tilde{u}_\lambda |_U$ to the map $u_0 |_{\overline{B}^2 \setminus U}$ with the help of 
an annulus embedded in a neighbourhood of $u_0 (\partial U)$. The obtained map can be made 
$J_\lambda$-holomorphic for some $J_\lambda$ which equals $J_0$ outside a neighbourhood of 
$u_0 (\partial U)$. We have then $\frac{d}{d\lambda} j_2 (u_\lambda , z_c) = 0$ and
$\frac{d}{d\lambda} j_3 (u_\lambda , z_c) = w$. Hence, $d_{(\stackrel{.}{\tilde{u}}_\lambda ,
\dot{J}_\lambda , z_c)} F_{deg} = \det (j_2 (z_c) , w)$ does not vanish as soon as
$w$ is not chosen colinear to $j_2 (z_c)$. $\square$

\begin{rem}
\label{remjet}
1) The group of real biholomorphisms of $(\overline{B}^2 , i)$ acts on $\R {\cal P}'_{cusp} (\eta)$
by $\phi .(u,J,z_c) = (u \circ \phi^{-1} , J , \phi (z_c))$. Since it is transitive on
$]-1,1[$, we can always assume that $z_c = 0$.

2) Let $(u,J, 0) \in \R {\cal P}'_{dcusp} (\eta)$ and $a \in \R$ be such that $j_3 (u,0) =
a j_2 (u,0)$. Let $z = w - \frac{a}{2}w^2$, the order three jet of $u$ writes
\begin{eqnarray*}
u(z) & =& j_2 (u,0) (w - \frac{a}{2}w^2)^2 + j_3 (u,0) (w - \frac{a}{2}w^2)^3 + o(|w|^3)\\
& =& j_2 (u,0) w^2 + o(|w|^3).
\end{eqnarray*}
\end{rem}

\begin{lemma}
\label{lemmajet}
Let $(u,J, 0) \in \R {\cal P}'_{dcusp} (\eta)$ and $a \in \R$ be such that $j_3 (u,0) =
a j_2 (u,0)$. Let $z = w - \frac{a}{2}w^2$, the order five germ of $u$ writes
$u(z) = j_2 w^2 + j_4^w w^4 + j_4^{|w|} |w|^4 + j_4^{\overline{w}} \overline{w}^4 +
j_5 w^5 + o(|w|^5),$ where $j_2 , j_4^w , j_4^{|w|} , j_4^{\overline{w}} , j_5 \in \R^2$.
\end{lemma}

{\bf Proof:}

This is deduced after expanding the relation $0= J|_u \circ du - du \circ i$. $\square$\\

A degenerated cuspidal curve $(u,J, z_c) \in \R {\cal P}'_{dcusp} (\eta)$ as in Lemma \ref{lemmajet}
is said to be {\it generic} if
the vectors $j_2$ and $j_5$ are linearly independant.

\begin{prop}
\label{proptransv}
Let $(u,J, 0) \in \R {\cal P}'_{dcusp} (\eta)$ be a generic cuspidal curve and 
$v(z) = d_z u(\frac{1}{z})$. Then
$(v,0,0) \in T_{(u,J, 0)} \R {\cal P}'_{cusp} (\eta) \setminus T \R {\cal P}'_{dcusp} (\eta)$.
Moreover, if $(u_\lambda, J_\lambda, z_c^\lambda)_{\lambda \in ]-\epsilon , \epsilon[}$ is a path of
$\R {\cal P}'_{cusp} (\eta)$ transversal to $\R {\cal P}'_{dcusp}$ at $\lambda = 0$, then
for $\lambda < 0$ (resp. $\lambda > 0$) close enough to $0$, $u_\lambda$ has an
isolated (resp. non isolated) real double point in a neighbourhood of the cusp
$u_\lambda (z_c^\lambda)$, or vice versa.
\end{prop}

{\bf Proof:}

With the notations of Lemma \ref{lemmajet}, we make the local change of coordinates
$z = \phi (w)$ with $\phi (w) = w - \frac{a}{2}w^2$. The order five germ of $u$ writes
$u \circ \phi (w) = j_2 w^2 + j_4^w w^4 + j_4^{|w|} |w|^4 + j_4^{\overline{w}} \overline{w}^4 +
j_5 w^5 + o(|w|^5),$ where $j_2 , j_4^w , j_4^{|w|} , j_4^{\overline{w}} , j_5 \in \R^2$.
Equip $\C^2$ with
the complex structure $J(0)$ and the frame $(u_0 (0) , j_2 , j_5)$. After composition with the
local diffeomorphism of $\C^2$ tangent to the identity given by $\Phi (z_1 , z_2) = (z_1 , z_2) - 
j_4^w z_1^2 -
j_4^{|w|} |z_1|^2 -  j_4^{\overline{w}} \overline{z}_1^2$, we can assume that the jet of
$u$ writes $j_2 w^2 + j_5 w^5 + o(|w|^5)$.
It suffices to prove that $\tilde{v} = d_w (\Phi \circ u \circ \phi) (\frac{1}{w}) \in 
T_{(u \circ \phi,J, 0)} \R {\cal P}'_{cusp} (\eta) \setminus T \R {\cal P}'_{dcusp} (\eta)$.
Indeed, $d_w (u \circ \phi) (\frac{1}{w}) = d_z u (\frac{1}{z} + b(z))$ for some
holomorphic $b$, and $d_z u (b(z)) \in T \R {\cal P}'_{dcusp} (\eta)$. Moreover, the composition by
$\Phi$ does not affect the transversality condition. Now
$\tilde{v} (w) = 2j_2 + 
5j_5 w^3 + o(|w|^3)$. From Lemma \ref{lemmarp'}, 
$(\tilde{v},0,0) \in T_{(u \circ \phi,J, 0)} \R {\cal P}'_{cusp} (\eta)$ and from Lemma
\ref{lemmarp'd}, $(\tilde{v},0,0) \notin T_{(u \circ \phi,J, 0)} \R {\cal P}'_{dcusp} (\eta)$,
since $dj_3 (\tilde{v},0,0) = 5j_5$ and $\det (j_2 , 5j_5) \neq 0$. The first part of the
proposition is proved.

Now, without loss of generality, we can assume that $\frac{d}{d\lambda} 
(u_\lambda, J_\lambda, z_c^\lambda) = (v,0,0)$ and $z_c^\lambda \equiv 0$. 
From what precedes, the expansion of $u_\lambda ( w - \frac{a}{2}w^2)$ writes
\begin{eqnarray*}
u_\lambda ( w - \frac{a}{2}w^2) &=& f(\lambda) + (w^2 , w^5) + O(|w|^6) \\
&&+ \lambda (0, 5 w^3 )+ O(\lambda |w|^4)\\
&&+ O(\lambda^2 |w|^2)
\end{eqnarray*}
The function $f(\lambda) + (w^2 , 5\lambda w^3 + w^5)$ has an isolated real double point at
the parameters $w = \pm i \sqrt{5\lambda}$ when $\lambda >0$. Let us prove that when
$\lambda >0$ is close enough to $0$, the function $u_\lambda ( w - \frac{a}{2}w^2)$ also
has an isolated real double point at parameters close to $\pm i \sqrt{5\lambda}$. Set
$w = i\sqrt{5\lambda} + \tilde{w}$, we have:
\begin{eqnarray*}
u_\lambda \circ \phi(i\sqrt{5\lambda} + \tilde{w}) & = & f(\lambda) + (-5\lambda , 0) +
\big( 2i\sqrt{5\lambda} \tilde{w} + \tilde{w}^2 + \lambda^2 O(|\tilde{w}|+|\lambda|) ,
50\lambda^2 \tilde{w} \\ 
&&+ |\lambda|^{\frac{3}{2}} \tilde{w}^2 O(|\tilde{w}|+|\lambda|)
 + \lambda^3 O(|\tilde{w}|+|\lambda|) \big),
\end{eqnarray*}
as soon as $|\tilde{w}| \ll \sqrt{\lambda}$. Set $\tilde{w} = |\lambda|^{\frac{3}{4}} 
(\cos(\theta) + i\sin (\theta))$, we get
$Im (u_\lambda \circ \phi(i\sqrt{5\lambda} + \tilde{w})) = \big( 2\sqrt{5}|\lambda|^{\frac{5}{4}}
\cos(\theta) + |\lambda|^{\frac{3}{2}} O(|\lambda|) , 50 |\lambda|^{\frac{11}{4}} \sin (\theta)
+ |\lambda|^3 O(|\lambda|) \big)$. The linking number between the origin of $\R^2$ and the ellipse
parameterized by $\theta \mapsto (2\sqrt{5}|\lambda|^{\frac{5}{4}}
\cos(\theta) , 50 |\lambda|^{\frac{11}{4}} \sin (\theta))$ is equal to one. The same result holds
for the linking number between the origin of $\R^2$ and the curve parameterized by 
$\theta \mapsto Im \big( u_\lambda \circ \phi(i\sqrt{5\lambda} + |\lambda|^{\frac{3}{4}} 
(\cos(\theta) + i\sin (\theta))) \big)$ as soon as $\lambda >0$ is close enough to $0$.
Hence, $Im (u_\lambda \circ \phi)$ vanishes once in the disk centered at $i\sqrt{5\lambda}$
and whose radius is $|\lambda|^{\frac{3}{4}}$. It follows that $u_\lambda$ has
an isolated real double point close to its cuspidal point $u_\lambda (0)$ as soon as
$\lambda >0$ is close enough to $0$.

In the same way, the function $f(\lambda) + (w^2 , 5\lambda w^3 + w^5)$ has a non isolated real 
double point at the parameters $w = \pm \sqrt{-5\lambda}$ when $\lambda <0$. Let us prove that when
$\lambda <0$ is close enough to $0$, the function $u_\lambda ( w - \frac{a}{2}w^2)$ also
has a non isolated real double point at parameters close to $\pm \sqrt{-5\lambda}$. Set
$w = \eta \sqrt{-5\lambda} + \tilde{w}$, where $\eta = \pm 1$, we have:
\begin{eqnarray*}
u_\lambda \circ \phi(\eta \sqrt{-5\lambda} + \tilde{w}) & = & f(\lambda) + (-5\lambda , 0) +
\big( 2\eta \sqrt{-5\lambda} \tilde{w} + \tilde{w}^2 + |\lambda|^2 O(|\tilde{w}|+|\lambda|) ,
50\lambda^2 \tilde{w} \\ 
&&+ |\lambda|^{\frac{3}{2}} \tilde{w}^2 O(|\tilde{w}|+|\lambda|)
 + \lambda^3 O(|\tilde{w}|+|\lambda|) \big),
\end{eqnarray*}
as soon as $|\tilde{w}| \ll \sqrt{|\lambda|}$. When $\tilde{w} \in [-|\lambda|^{\frac{3}{4}} ,
|\lambda|^{\frac{3}{4}} ]$, the segment
$(2\eta \sqrt{-5\lambda} \tilde{w} , 50\lambda^2 \tilde{w})$ joins the two points
$(-2\sqrt{5} \eta |\lambda|^{\frac{5}{4}} , -50 |\lambda|^{\frac{11}{4}})$ and
$(2\sqrt{5} \eta |\lambda|^{\frac{5}{4}} , 50 |\lambda|^{\frac{11}{4}})$. When $\eta = \pm 1$,
these two segments intersect transversely at the origin. The same is true for the segments
$u_\lambda \circ \phi(\eta \sqrt{-5\lambda} +[-|\lambda|^{\frac{3}{4}} ,
|\lambda|^{\frac{3}{4}} ])$ as soon as $\lambda <0$ is close enough to $0$, hence the
result. $\square$

\subsection{Global study of degenerated order two cuspidal points}
\label{subsecglobalstudy}

Denote by $\R {\cal M}^d_{dcusp}$ the subset of elements $[u , J_S , J , \underline{x} , z_c]$
of $\R {\cal M}^d_{cusp}$ for which $u$ has a degenerated cuspidal point at $z_c$.
\begin{prop}
\label{propdegcusp}
The space $\R {\cal M}^d_{dcusp}$ is a codimension one submanifold of $\R {\cal M}^d_{cusp}$
of class $C^{l-k}$. Moreover, the subspace of $\R {\cal M}^d_{dcusp}$ made of curves
$[u , J_S , J , \underline{x} , z_c]$ which have a non generic degenerated cuspidal point at $z_c$
is a substratum of codimension two and class $C^{l-k}$ of $\R {\cal M}^d_{cusp}$. $\square$
\end{prop}
(This is a particular case of Theorem $3.4.5$ of \cite{Shev}, see also Proposition $7$ of \cite{Bar})\\

Let $\gamma : [0,1] \to \R {\cal J}_\omega \times \R_\tau X^{c_1(X)d - 2}$ be a path transversal 
to $\pi_\R$. From Proposition \ref{propdegcusp}, if it is chosen generic enough, it
avoids the image under $\pi_\R$ of curves having a non generic degenerated cuspidal point.
Denote by $\R {\cal M}_\gamma = \R {\cal M}^d_{cusp} \times_\gamma [0,1]$ and by
$\pi_\gamma : \R {\cal M}_\gamma \to [0,1]$ the associated projection.
\begin{prop}
\label{propcritdegcusp}
Let $C_{t_0} = [u , J_S , J^{t_0} , \underline{x}^{t_0} , z_c] \in \R {\cal M}_\gamma$ be a
curve having a degenerated cuspidal point at $z_c$. Then, there exists a neighbourhood $W$
of $[u , J_S , J^{t_0} , \underline{x}^{t_0} , z_c]$ in $\R {\cal M}_\gamma$ and $\epsilon > 0$
such that for every $t \in ]t_0 - \epsilon , t_0[$, $\pi_\gamma^{-1} (t) \cap W$ is made of
two curves $C_t^+$, $C_t^-$ such that $m(C_t^+) = m(C_t^-) + 1$ and for every 
$t \in ]t_0 , t_0 + \epsilon [$, $\pi_\gamma^{-1} (t) \cap W = \emptyset$, or vice versa.
\end{prop}
{\bf Proof:}

From Proposition \ref{propdegcusp}, the degenerated cuspidal point of $C_{t_0}$ is generic. 
From Theorem \ref{theocrit}, $C_{t_0}$ is a non-degenerated critical point of $\pi_\gamma$. Thus,
as soon as $\epsilon$ and $W$ are small enough, for every $t \in ]t_0 - \epsilon , t_0[$, 
$\pi_\gamma^{-1} (t) \cap W = \{ C_t^\pm \}$ and for every 
$t \in ]t_0 , t_0 + \epsilon [$, $\pi_\gamma^{-1} (t) \cap W = \emptyset$, or vice versa. The
only thing to prove is that $m(C_t^+) = m(C_t^-) + 1$. The double points of $C_t^\pm$ are
close to the ones of $C_{t_0}$ with the exception of one which is close to the cusp of
$C_{t_0}$. We have to prove that the nature of the latter is not the same for
$C_t^+$ and $C_t^-$. Note that $\frac{d}{dt}|_{t=t_0} C_t$ is the generator of $\ker d\pi_\R =
H^0 (S , {\cal N}^{cusp}_{u , -\underline{z}})_{+1} \cong H^0 (S , \C_{z_c , -\underline{z}})_{+1}$, see
Proposition \ref{propprepcrit}. Choose a neighbourhood $V$ of $u(z_c)$ invariant under $c_X$, diffeomorphic
to the $4$-ball $B^4$ and small enough in order that $||J^t - J^{t_0} (u (z_c))||_{C^1} < \eta$.
We deduce a restriction map $rest : W \subset \R {\cal M}_\gamma \to \R {\cal P}'_{cusp} (\eta)$
such that $rest (C_{t_0}) \in \R {\cal P}'_{dcusp} (\eta)$. Now 
$d_{C_{t_0}} rest (\frac{d}{dt}|_{t=t_0} C_t)$ is exactly the vector $(v,0,0)$ given by
Proposition \ref{proptransv}. The result thus follows from Proposition \ref{proptransv}.
$\square$

\section{Study of the compactification $\R \overline{\cal M}^d_{cusp}$}
\label{sectioncompact}

\subsection{Neighbourhood of $\R {\cal M}^{(d_1 , m_1), (d_2 , m_2)}_{tac}$ in 
$\R \overline{\cal M}^d_{cusp}$ when $m_i = c_1 (X) d_i - 1$}
\label{subsecttac1}

\begin{lemma}
\label{lemmacoorient}
Let $d_1$, $d_2 \in H_2 (X ; \Z)$ be such that $d_1 + d_2 = d$ and $m_i = c_1 (X) d_i -1$,
$i \in \{ 1, 2 \}$.
Then, $\R {\cal M}^{(d_1 , m_1), (d_2 , m_2)}_{tac}$ has a canonical coorientation in
$\R {\cal M}^{(d_1 , m_1), (d_2 , m_2)}$. A path $(C^t)_{t \in ]-\epsilon, \epsilon[}$ in
$\R {\cal M}^{(d_1 , m_1), (d_2 , m_2)}$ transversal to 
$\R {\cal M}^{(d_1 , m_1), (d_2 , m_2)}_{tac}$ at $t = 0$ is positive on this coorientation
if for every $t \in ]-\epsilon, 0[$ (resp. $t \in ]0 , \epsilon [$), $C^t$ has two real
( resp. complex conjugated) double points in a neighbourhood of the tacnode of $C^0$.
$\square$
\end{lemma}
(See \S \ref{subsectcompact} for the definition of the space
$\R {\cal M}^{(d_1 , m_1), (d_2 , m_2)}_{tac}$)
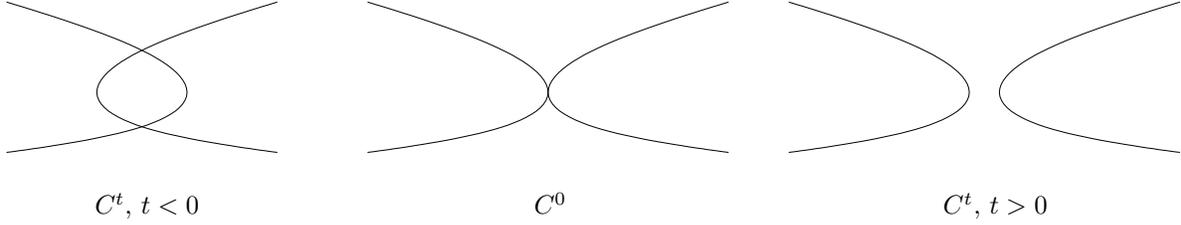
$$\vcenter{\hbox{\begin{picture}(0,0)%
\includegraphics{cusp1.pstex}%
\end{picture}%
\setlength{\unitlength}{3315sp}%
\begingroup\makeatletter\ifx\SetFigFont\undefined%
\gdef\SetFigFont#1#2#3#4#5{%
  \reset@font\fontsize{#1}{#2pt}%
  \fontfamily{#3}\fontseries{#4}\fontshape{#5}%
  \selectfont}%
\fi\endgroup%
\begin{picture}(8799,1587)(1564,-6811)
\put(5536,-6811){\makebox(0,0)[lb]{\smash{\SetFigFont{10}{12.0}{\rmdefault}{\mddefault}{\updefault}{$C^0$}%
}}}
\put(8596,-6811){\makebox(0,0)[lb]{\smash{\SetFigFont{10}{12.0}{\rmdefault}{\mddefault}{\updefault}{$C^t, \, t>0$}%
}}}
\put(2251,-6811){\makebox(0,0)[lb]{\smash{\SetFigFont{10}{12.0}{\rmdefault}{\mddefault}{\updefault}{$C^t , \, t<0$}%
}}}
\end{picture}
}}$$
Remember that under the conditions of Lemma \ref{lemmacoorient}, the projection
$\pi_\R^{(d_1 , m_1), (d_2 , m_2)} : \R {\cal M}^{(d_1 , m_1), (d_2 , m_2)} \to
\R {\cal J}_\omega \times \R_{\tau_1} X^{m_1} \times \R_{\tau_2} X^{m_2}$ is Fredholm
with vanishing index. Moreover, from Proposition \ref{propredrap}, 
$\R {\cal M}^{(d_1 , m_1), (d_2 , m_2)}_{tac}$ is made of regular points of this projection.
Let $\gamma : t \in [0,1] \mapsto (J^t , \underline{x}^t) \in \R {\cal J}_\omega \times 
\R_{\tau_1} X^{m_1} \times \R_{\tau_2} X^{m_2}$ be a path transversal to 
$\pi_\R^{(d_1 , m_1), (d_2 , m_2)}$, so that the fibre product $\R {\cal M}^{red}_\gamma =
\R {\cal M}^{(d_1 , m_1), (d_2 , m_2)} \times_{\pi_\R} [0,1]$ is a smooth one dimensional
manifold equipped with a projection $\pi^{red}_\gamma : \R {\cal M}^{red}_\gamma \to [0,1]$.
As soon as $\gamma$ is chosen generic enough, this submanifold intersects
$\R {\cal M}^{(d_1 , m_1), (d_2 , m_2)}_{tac}$ transversely at finitely many points over distinct
parameters of $[0,1]$. Let $C^{t_0} \in  \R {\cal M}^{red}_\gamma$ be such a point,
$t_0 \in ]0,1[$. The path $\gamma$ is said to intersect $\pi_\R^{(d_1 , m_1), (d_2 , m_2)} 
(\R {\cal M}^{(d_1 , m_1), (d_2 , m_2)}_{tac})$ {\it positively} (resp. {\it negatively})
at $\gamma (t_0)$ if $\R {\cal M}^{red}_\gamma$ intersects 
$\R {\cal M}^{(d_1 , m_1), (d_2 , m_2)}_{tac}$ positively (resp. negatively) at
$C^{t_0}$ once endowed with the local orientation induced by $[0,1]$.

Assume that $\gamma$ is transversal to $\pi_\R$ and denote by $\R {\cal M}_\gamma =
\R {\cal M}^d_{cusp} \times_{\gamma} [0,1]$. Denote by $\R \overline{\cal M}_\gamma$ the Gromov
compactification of $\R {\cal M}_\gamma$ and by $\pi_\gamma : \R \overline{\cal M}_\gamma \to
[0,1]$ the associated projection. The aim of this subparagraph is to prove the following
theorem.
\begin{theo}
\label{theoredtac1}
Let $\gamma : t \in [0,1] \mapsto (J^t , \underline{x}^t) \in \R {\cal J}_\omega \times 
\R_{\tau} X^{c_1 (X)d - 2}$ be a generic path chosen as above and
$C^{t_0} \in \R \overline{\cal M}_\gamma \cap \R {\cal M}^{(d_1 , m_1), (d_2 , m_2)}_{tac}$.
Assume that $m_1 = c_1 (X) d_1 -1$, $m_2 = c_1 (X) d_2 -1$ and that $\gamma$ is positive at
$t_0 = \pi_\gamma (C^{t_0})$. Then, there exist a neighbourhood $W$ of $C^{t_0}$ in 
$\R \overline{\cal M}_\gamma$ and $\epsilon > 0$ such that for every 
$t \in ]t_0 - \epsilon , t_0[$, $\pi_\gamma^{-1} (t) \cap W$ is made of two curves having the
same mass, the one of $C^{t_0}$, and for every 
$t \in ]t_0 , t_0 + \epsilon[$, $\pi_\gamma^{-1} (t) \cap W = \emptyset$.
\end{theo}
Note that reversing the orientation of $[0,1]$ if necessary, we can always assume that
$\gamma$ is positive at $t_0$.

Let $C^t$ be a real rational cuspidal $J^t$-holomorphic curve close to $C^{t_0}$ which passes
through $\underline{x}^t$, $t \in ]t_0 - \epsilon , t_0 + \epsilon[ \setminus \{ t_0 \}$. Then,
from Proposition $2.16$ of \cite{Wels1}, $C^t$ extends to a one parameter family of 
$J^t$-holomorphic curves $C^t (\eta)$, $\eta \in [0 , \eta_t [$ such that $C^t (0) =
C^t$, $C^t (\eta)$ passes
through $\underline{x}^t$ for every $\eta \in ]0 , \eta_t [$ and $\R C^t (\eta)$ has a
non isolated real double point in the neighbourhood of the cusp of $C^t$ as soon as
$\eta \neq 0$.
$$\vcenter
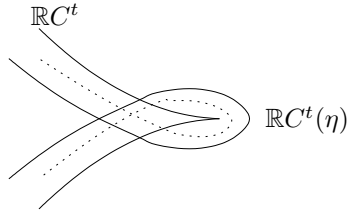
{\hbox{\begin{picture}(0,0)%
\includegraphics{cusp2.pstex}%
\end{picture}%
\setlength{\unitlength}{3315sp}%
\begingroup\makeatletter\ifx\SetFigFont\undefined%
\gdef\SetFigFont#1#2#3#4#5{%
  \reset@font\fontsize{#1}{#2pt}%
  \fontfamily{#3}\fontseries{#4}\fontshape{#5}%
  \selectfont}%
\fi\endgroup%
\begin{picture}(1947,1572)(2239,-5473)
\put(4186,-4831){\makebox(0,0)[lb]{\smash{\SetFigFont{10}{12.0}{\rmdefault}{\mddefault}{\updefault}{$\R C^t (\eta)$}%
}}}
\put(2431,-4066){\makebox(0,0)[lb]{\smash{\SetFigFont{10}{12.0}{\rmdefault}{\mddefault}{\updefault}{$\R C^t$}%
}}}
\end{picture}
}}$$
\begin{lemma}
\label{lemmafamily}
If $\epsilon$ is small enough, the family $C^t (\eta)$ converges to a reducible
$J^t$-holomorphic curve when $\eta \to \eta_t$.
\end{lemma}
Note that Lemma \ref{lemmafamily} already implies the second part of Theorem \ref{theoredtac1},
since for $t > 0$, there are no reducible
$J^t$-holomorphic curve which pass
through $\underline{x}^t$ and have a real double point in the neighbourhood of the tacnode
 of $C^{t_0}$.\\

{\bf Proof:}

Let $U$ be a compact neighbourhood of $C^{t_0}$ in $X$ such that for every 
$t \in ]t_0 - \epsilon , t_0 + \epsilon[$, the only reducible $J^t$-holomorphic curve which pass
through $\underline{x}^t$ are the ones close to $C^{t_0}$. Note that as soon as $\eta$ is
close enough to zero, the real parts $\R C^t (\eta)$ form a loop around the cusp of $\R C^t$
in $\R X$. Moreover, the intersections between two curves of this family $C^t (\eta)$ are
located at $\underline{x}^t$ and in the neighbourhood of their double points. Thus, as $\eta$
grows, the loops grow in order to fill a disk of $\R X$ centered at the cusp of $\R C^t$. The
following alternative now
follows from Gromov's compactness Theorem. Either $C^t (\eta)$ converges to a reducible
$J^t$-holomorphic curve in $U$ as $\eta \to \eta_t$, or one curve $C^t (\eta)$ intersects the
boundary of $U$. Assume that there exists a sequence $t_n \in 
]t_0 - \epsilon , t_0 + \epsilon[ \setminus \{ t_0 \}$, $n \in \N^*$, which converges to $t_0$ when 
$n \to \infty$ and $\eta_n >0$ such that $C^{t_n} (\eta_n)$ touches the
boundary of $U$. Then, when $n \to \infty$, $C^{t_n} (\eta_n)$ converges to a 
$J^{t_0}$-holomorphic curve $C^\infty$ which is contained in $U$, intersects the
boundary of $U$ and passes through $\underline{x}^{t_0}$. The latter cannot be reducible from 
the definition of $U$. Moreover, for every $n \in \N^*$, the loop of $\R C^{t_n} (\eta_n)$
surounds the cusp of $\R C^{t_n}$ in $\R X$. It follows that in the limit, $\R C^\infty$ forms
a loop which surounds the tacnode of $\R C^{t_0}$. Thus, $C^\infty$ intersects $C^{t_0}$ with
multiplicity four at least near the tacnode of $C^{t_0}$, with multiplicity two at least near 
every double point of $C^{t_0}$ and with multiplicity one at $\underline{x}^{t_0}$. The total
intersection index between $C^\infty$ and $C^{t_0}$ is then at least $d^2 +2$, which is
impossible. $\square$
\begin{lemma}
\label{lemmatwo}
Assume that $\epsilon$ is small enough and that $t \in ]t_0 - \epsilon , t_0 [$. Then, the number
of cuspidal $J^t$-holomorphic curves which pass through $\underline{x}^t$ and are close
to $C^{t_0}$ is at most $2$.
\end{lemma}
{\bf Proof:}

Denote by $C^t_{red}$ the unique reducible real rational $J^t$-holomorphic curve which passes 
through $\underline{x}^t$ and is close to $C^{t_0}$. This curve has two non isolated real
double points $y_1^t$, $y_2^t$ in a neighbourhood of the tacnode of $C^{t_0}$. Let $C^t_1$ be 
a real rational cuspidal $J^t$-holomorphic curve which passes 
through $\underline{x}^t$ and is close to $C^{t_0}$ and $C^t_1 (\eta)$, $\eta \in [0,\eta_1^t]$,
the one parameter family of $J^t$-holomorphic curves given by Lemma \ref{lemmafamily}. In
particular, $C^t_1 (\eta_1^t) = C^t_{red}$. For every $\eta \in ]0,\eta_1^t[$, denote by
$y_1^t (\eta)$ the real double point of $C^t_1 (\eta)$ close to the tacnode of $C^{t_0}$.
The latter converges to the cusp of $C^t_1$ when $\eta \to 0$ and to one of the points
$y_1^t$, $y_2^t$, say $y_1^t$, when $\eta \to \eta_1^t$. Assume that there were two families
$C^t_1 (\eta)$, $C^t_2 (\delta)$ of curves having this property. Then, for $\eta, \delta$
close to $\eta_1^t$, $\delta_1^t$, the curves $C^t_1 (\eta)$, $C^t_2 (\delta)$ would have all
their intersections at $\underline{x}^t$ and in the neighbourhood of the double points of
$C^t_{red}$ different from $y_2^t$. Moreover, if $\eta$ is close enough to $\eta_1^t$, we can
assume that the loop formed by $\R \C^t_2 (\delta)$ close to the cusp of $C^t_1$ is in the
interior of the one formed by $\R C^t_1 (\eta)$.
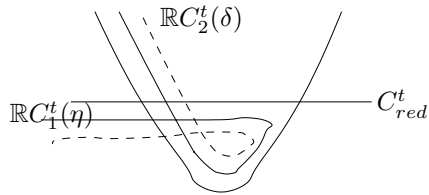
$$\vcenter{\hbox{\begin{picture}(0,0)%
\includegraphics{cusp3.pstex}%
\end{picture}%
\setlength{\unitlength}{3315sp}%
\begingroup\makeatletter\ifx\SetFigFont\undefined%
\gdef\SetFigFont#1#2#3#4#5{%
  \reset@font\fontsize{#1}{#2pt}%
  \fontfamily{#3}\fontseries{#4}\fontshape{#5}%
  \selectfont}%
\fi\endgroup%
\begin{picture}(2745,1467)(1576,-4798)
\put(2701,-3526){\makebox(0,0)[lb]{\smash{\SetFigFont{10}{12.0}{\rmdefault}{\mddefault}{\updefault}{$\R C^t_2 (\delta)$}%
}}}
\put(1576,-4246){\makebox(0,0)[lb]{\smash{\SetFigFont{10}{12.0}{\rmdefault}{\mddefault}{\updefault}{$\R C^t_1 (\eta)$}%
}}}
\put(4321,-4156){\makebox(0,0)[lb]{\smash{\SetFigFont{10}{12.0}{\rmdefault}{\mddefault}{\updefault}{$C^t_{red}$}%
}}}
\end{picture}
}}$$
Then, $\R C^t_1 (\eta)$ intersects $\R \C^t_2 (\delta)$ at two points belonging to the two
local branches of $\R \C^t_2 (\delta)$ near $y_1^t$. As $\eta$ decreases, there is
some parameter $\eta'$ for which $\R C^t_1 (\eta')$ passes through the double point of
$\R \C^t_2 (\delta)$ belonging to its loop. For this parameter $\eta'$, the two curves
$C^t_1 (\eta')$ and $\C^t_2 (\delta)$ would have at least $3$ intersection points in the 
neighbourhood of the tacnode of $C^{t_0}$ and thus a total intersection at least equal to
$d^2 + 1$, which is impossible. We deduce that the number of real rational cuspidal 
$J^t$-holomorphic curves close to $C^{t_0}$ is bounded by the number of real double points
of $C^t_{red}$ close to the tacnode of $C^{t_0}$, that is two. $\square$\\

{\bf Proof of the Theorem \ref{theoredtac1}:}

It follows from Lemmas \ref{lemmafamily} and \ref{lemmatwo} that there exist
a neighbourhood $W$ of $C^{t_0}$ in 
$\R \overline{\cal M}_\gamma$ and $\epsilon > 0$ such that for every 
$t \in ]t_0 - \epsilon , t_0[$, $\# (\pi_\gamma^{-1} (t) \cap W) \leq 2$ and for every 
$t \in ]t_0 , t_0 + \epsilon[$, $\pi_\gamma^{-1} (t) \cap W = \emptyset$. Since the parity of
the number of real rational cuspidal $J$-holomorphic curves which pass through $\underline{x}$
does not depend on the generic choice of $(J,\underline{x})$, it suffices to prove that
$\pi_\gamma^{-1} (t) \cap W$ cannot be empty when $t \in ]t_0 - \epsilon , t_0[$, that is there
exists at least one real rational cuspidal $J$-holomorphic curves which passes through 
$\underline{x}$ for some $J$ close to $J^{t_0}$. Now such a curve can be constructed
by reversing the construction of Lemma \ref{lemmafamily}. One starts with a reducible 
real $J$-holomorphic curve having two real non isolated double points close to the
tacnode of $C^{t_0}$. Then, from Proposition $2.14$ of \cite{Wels1}, 
one can smooth one of these double points
to obtain a one parameter family $\R C (\eta)$ of curves which forms a loop. As $\eta$
decreases, it has to degenerate onto a cuspidal curve by some argument similar to the
one used in the proof of Lemma \ref{lemmafamily}. $\square$

\begin{rem}
\label{rem0}
When $r = (0, \dots , 0)$, a sequence of real rational cuspidal
$J$-holomorphic curves can converge to a reducible curve $C$ having two irreducible components which
are complex conjugated and transversal to each other except at one point which is of order two, that
is an ordinary tacnode. To extend Theorem \ref{maintheo}, one should take into account these
reducible curves. This would be possible provided an analog of Theorem \ref{theoredtac1} holds in this case.
Namely, assume that over a path $\gamma$, the curve deforms to a reducible curve having
two real (resp. complex conjugated) double points in a neighbourhood of the tacnode for
$t \in ] t_0 - \epsilon , t_0 [$ (resp. $t \in ]t_0 , t_0 + \epsilon[$). Then, one can suspect that
for $t \in ]t_0 , t_0 + \epsilon[$, there are no real rational cuspidal $J^t$-holomorphic curve close
to $C$ whereas there are two of them for $t \in ] t_0 - \epsilon , t_0 [$. Moreover, the latter come from
the two degenerations of the figure eight. However, I have no proof of this fact and thus leave this case 
opened.
\end{rem}

\subsection{Neighbourhood of $\R {\cal M}^{(d_1 , m_1), (d_2 , m_2)}_{tac}$ in 
$\R \overline{\cal M}^d_{cusp}$ when $m_2 = c_1 (X) d_2$}
\label{subsecttac2}

Let $\gamma : t \in [0,1] \mapsto (J^t , \underline{x}^t) \in \R {\cal J}_\omega \times 
\R_{\tau} X^{c_1 (X)d - 2}$ be a generic path 
transversal to $\pi_\R$ and $\R {\cal M}_\gamma =
\R {\cal M}^d_{cusp} \times_{\gamma} [0,1]$. Denote by $\R \overline{\cal M}_\gamma$ the Gromov
compactification of $\R {\cal M}_\gamma$ and by $\pi_\gamma : \R \overline{\cal M}_\gamma \to
[0,1]$ the associated projection. 
\begin{theo}
\label{theoredtac2}
Let $\gamma : t \in [0,1] \mapsto (J^t , \underline{x}^t) \in \R {\cal J}_\omega \times 
\R_{\tau} X^{c_1 (X)d - 2}$ be a generic path chosen as above and
$C^{t_0} \in \R \overline{\cal M}_\gamma \cap \R {\cal M}^{(d_1 , m_1), (d_2 , m_2)}_{tac}$.
Assume that $m_1 = c_1 (X) d_1 -2$ and $m_2 = c_1 (X) d_2$. Then, there exist a neighbourhood $W$ of 
$C^{t_0}$ in 
$\R \overline{\cal M}_\gamma$ and $\epsilon > 0$ such that for every 
$t \in ]t_0 - \epsilon , t_0 + \epsilon[$, $\pi_\gamma^{-1} (t) \cap W = \{ C^t \}$. 
Moreover, the mass of $C^t$ does not depend on $t \in ]t_0 - \epsilon , t_0 + \epsilon[$.
\end{theo}
(See \S \ref{subsectcompact} for the definition of the space
$\R {\cal M}^{(d_1 , m_1), (d_2 , m_2)}_{tac}$)

Without loss of generality, we can assume that $\underline{x}^t \in \R_{\tau} X^{c_1 (X)d - 2}$
does not depend on $t \in [0,1]$. Denote by $C^{t_0}_1$ (resp $C^{t_0}_2$) the irreducible
component of $C^{t_0}$ which has homology class $d_1$ (resp. $d_2$). Let $\underline{x}^t_1 =
\underline{x}^t \cap C^{t_0}_1$, $\underline{x}^t_2 = \underline{x}^t \cap C^{t_0}_2$ and
$U$ be a compact neighbourhood of $C^{t_0}$ in $X$. If $U$ is small enough, the curve $C^{t_0}_1$
extends to a one parameter family $C^{t_0}_1 (\eta)$, $\eta \in [-1 , 1]$, of real
$J^{t_0}$-holomorphic curves which pass through $\underline{x}^t_1$, are contained in $U$ and
such that $C^{t_0}_1 (0) = C^{t_0}_1$, $C^{t_0}_1 (\pm 1) \cap \partial U \neq \emptyset$
whereas $C^{t_0}_1 (\eta) \subset \stackrel{\circ}{U}$ for $\eta \in ]-1 , 1[$. We
can assume that the curves $\R C^{t_0}_1 (\eta)$ and $\R C^{t_0}_2$ have two real intersection 
points in a neighbourhood of the tacnode of $C^{t_0}$ when $\eta > 0$. Let $C^t$ be a
real rational cuspidal $J^t$-holomorphic curve close to $C^{t_0}$ and which passes 
through $\underline{x}^t$, $t \in ]t_0 - \epsilon , t_0 + \epsilon[ \setminus \{ t_0 \}$. Then,
from Proposition $2.16$ of \cite{Wels1}, $C^t$ extends to a one parameter family of 
$J^t$-holomorphic curves $C^t (\eta)$, $\eta \in [0 , 1]$ such that $C^t (0) =
C^t$, $C^t (\eta)$ passes
through $\underline{x}^t$ for every $\eta \in [0 , 1]$ and $\R C^t (\eta)$ has a
non isolated real double point in the neighbourhood of the cusp of $C^t$ as soon as
$\eta \neq 0$. Note that in contrast with \S \ref{subsecttac1}, as soon as $U$ is small enough,
this family cannot break into a reducible curve as long as it stays in $U$. We can thus assume that 
$C^t (1) \cap \partial U \neq \emptyset$.

\begin{lemma}
\label{lemmafamily2}
As soon as $U$ is small enough, $C^t (1)$ converges to $C^{t_0} (1)$ as $t$ converges to $t_0$.
\end{lemma}

{\bf Proof:}

It suffices to prove that as $t$ converges to $t_0$, the curve $C^t (1)$ converges to a reducible
curve. Indeed, since this curve in the limit is contained in $U$, touches $\partial U$ and
passes through $\underline{x}^{t_0}$, it has to coincide with $C^{t_0} (1)$. Remember
that when $t \in ]t_0 - \epsilon , t_0 + \epsilon[ \setminus \{ t_0 \}$, the curve
$\R C^t (1)$ forms a loop which surrounds the cusp of $\R C^t$. Hence, if the curve in the limit
were irreducible, its real part would form a loop which would surround the tacnode of
$C^{t_0}$. As soon as $U$ is small enough, the latter would then intersect $C^{t_0}$ with
multiplicity four near the tacnode, with
multiplicity two near every double point of $C^{t_0}$ and at $\underline{x}^{t_0}$, which is
impossible. $\square$\\

Hence, the family $C^t (\eta)$, $\eta \in ]0,1]$, is obtained after smoothing one of the two
real double points of $C^{t_0} (\eta)$ close to the tacnode of $C^{t_0}$.
\begin{lemma}
\label{lemmaone}
Let $C^t (\eta)$ and $C^t (\eta)'$, $t \in ]t_0 - \epsilon , t_0 + \epsilon[ \setminus \{ t_0 \}$,
$\eta \in [0,1]$, be two families of real rational $J^t$-holomorphic curves which pass
through $\underline{x}^t$ and such that $C^t (0)$, $C^t (0)'$ are cuspidal. Then, these families
are obtained after smoothing the same real double point of the family $C^{t_0} (\eta)$.
\end{lemma}

{\bf Proof:}

Assume the converse. Then, the curves $C^t (1)$ and $C^t (\frac{1}{2})'$ would intersect at
$\underline{x}^t$ and with multiplicity two near every double point of $C^{t_0}$. Moreover, the
loops formed by $\R C^t (1)$ and $\R C^t (\frac{1}{2})'$ would intersect at two points
at least. Finally, since $C^t (\frac{1}{2})'$ is not obtained after smoothing the same real 
double point of the family $C^{t_0} (\eta)$ as $C^t (1)$, the local real branch in this smoothing
which is not included in the loop of $\R C^t (\frac{1}{2})'$ would also intersect 
$\R C^t (1)$ near the tacnode of $C^{t_0}$. This would provide a total intersection index
at least equal to $c_1 (X)d - 2 + d^2 - c_1 (X)d + 3 > d^2$, hence the contradiction. $\square$
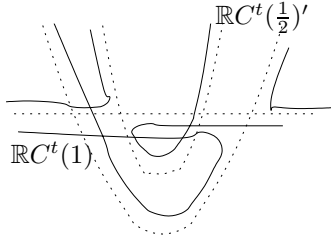
$$\vcenter{\hbox{\begin{picture}(0,0)%
\includegraphics{cusp4.pstex}%
\end{picture}%
\setlength{\unitlength}{3315sp}%
\begingroup\makeatletter\ifx\SetFigFont\undefined%
\gdef\SetFigFont#1#2#3#4#5{%
  \reset@font\fontsize{#1}{#2pt}%
  \fontfamily{#3}\fontseries{#4}\fontshape{#5}%
  \selectfont}%
\fi\endgroup%
\begin{picture}(2454,1782)(2194,-5248)
\put(2251,-4696){\makebox(0,0)[lb]{\smash{\SetFigFont{10}{12.0}{\rmdefault}{\mddefault}{\updefault}{$\R C^t (1)$}%
}}}
\put(3781,-3661){\makebox(0,0)[lb]{\smash{\SetFigFont{10}{12.0}{\rmdefault}{\mddefault}{\updefault}{$\R C^t (\frac{1}{2})'$}%
}}}
\end{picture}
}}$$

\begin{rem}
\label{remsame}
It also follows from the proof of Lemma \ref{lemmaone} that if the families $C^t (\eta)$ and 
$C^t (\eta)'$ are obtained after smoothing the same real double point of the family 
$C^{t_0} (\eta)$, $\eta \in [0,1]$, then the loop formed by $\R C^t (\frac{1}{2})'$ is
included in the one formed by $\R C^t (1)$. Indeed, these curves would otherwise also intersect
with  multiplicity at least three near the tacnode of $C^{t_0}$, which is impossible.
\end{rem}

{\bf Proof of Theorem \ref{theoredtac2}:}

There exist a neighbourhood $W$ of $C^{t_0}$ in 
$\R \overline{\cal M}_\gamma$ and $\epsilon > 0$ such that for every 
$t \in ]t_0 - \epsilon , t_0 + \epsilon[$, $\# (\pi_\gamma^{-1} (t) \cap W) \leq 1$. Indeed, if this
set would contain two curves, they would generate two families $C^t (\eta)$ and $C^t (\eta)'$
as in Lemmas \ref{lemmafamily2} and \ref{lemmaone}. From Lemma \ref{lemmaone}, these two families
would be obtained after smoothing the same real double point of the family $C^{t_0} (\eta)$.
From Remark \ref{remsame}, the loop formed by $\R C^t (\frac{1}{2})'$ and hence by
$\R C^t (\eta)'$ for every $\eta \in ]0,1]$ would be included in the one of $\R C^t (1)$. We
then obtain a contradiction repeating the proof of Lemma \ref{lemmatwo}. Moreover, as in the proof
of Theorem \ref{theoredtac1}, $\pi_\gamma^{-1} (t) \cap W$ cannot be empty for every
$t \in ]t_0 - \epsilon , t_0 + \epsilon[$, $\epsilon$ small enough. Since the parity of the
cardinality of $\pi_\gamma^{-1} (t) \cap W$ does not depend on $t$ and the masses of cuspidal
curves close to $C^{t_0}$ are obviously the ones of $C^{t_0}$, the theorem is proved. $\square$

\subsection{Neighbourhood of $\R {\cal M}^{(d_1 , m_1), (d_2 , m_2)}_{cusp}$ in 
$\R \overline{\cal M}^d_{cusp}$}
\label{subsectcusp}

\begin{prop}
\label{propinteg}
Let $(C^0 , J^0) \in \R {\cal M}^{(d_1 , m_1), (d_2 , m_2)}_{cusp}$. Then, there exists a
path $(C^t , J^t) \in \R {\cal M}^{(d_1 , m_1), (d_2 , m_2)}_{cusp}$, $t \in [0,1]$, of 
class $C^{l-k}$, such that $J^1$ is integrable in a neighbourhood of $C^1$ in $X$.
\end{prop}
(See \S \ref{subsectcompact} for the definition of the space
$\R {\cal M}^{(d_1 , m_1), (d_2 , m_2)}_{cusp}$)

\begin{lemma}
\label{lemmapath}
Under the hypothesis of Proposition \ref{propinteg}, there exists a neighbourhood $V$ of
the singular points of $C^0$ in $X$ and a
path $(C^t , J^t) \in \R {\cal M}^{(d_1 , m_1), (d_2 , m_2)}_{cusp}$, $t \in [0,1]$, of 
class $C^{l-k}$, such that $J^1|_V$ is integrable.
\end{lemma}

{\bf Proof:}

Let $y_0$ be a double point of $C^0$. Then, there exists a neighbourhood $V_0$ of $y_0$ in $X$,
invariant under $c_X$, as well as a diffeomorphism $\Phi : V_0 \to B^4$ which is 
$\Z/2\Z$-equivariant such that $\Phi (C^0 \cap V_0) = \{ (w_1 , w_2) \in B^4 \subset \C^2 \, | \,
w_1 w_2 = 0 \}$. Denote by $J$ the almost complex structure $\Phi_* (J^0|_{V_0})$ of $B^4$,
and for $\delta \in ]0,1]$, by $h_\delta$ the homothety $(w_1 , w_2) \in B^4 \mapsto
(\delta w_1 , \delta w_2) \in B^4$. Set then $J_{1 - \delta} = h_\delta^* (J)$, so that for
$\xi \in \R^4 = T_{(w_1 , w_2)} B^4$, $J_{1 - \delta}|_{(w_1 , w_2)} (\xi) =
dh_\delta^{-1} \circ J_{(\delta w_1 , \delta w_2)} \circ dh_\delta (\xi) = 
J_{(\delta w_1 , \delta w_2)} (\xi)$. Hence, the family $(J_{1 - \delta})_{\delta \in 
]0,1]}$, extend to a $C^\infty$-family $(J_\delta)_{\delta \in [0,1]}$, by setting
$J_1 \equiv J_{(0,0)}$. Now, $\Phi (C^0 \cap V_0)$ is invariant under $h_\delta$ and thus
$J_\delta$-holomorphic for every $\delta \in [0,1]$. The family 
$(\Phi^* J_\delta)_{\delta \in [0,1]}$ is then a path of almost complex structures of class $C^l$
of $V_0$ such that $\Phi^* J_0 = J^0|_{V_0}$, $C^0 \cap V_0$ is $\Phi^* J_\delta$-holomorphic
for every $\delta \in [0,1]$ and $\Phi^* J_1$ is integrable. There is no obstruction to extend
$(\Phi^* J_\delta)_{\delta \in [0,1]}$ to a path of almost complex structures of class $C^l$ on
the whole $X$ for which $C^0$ is holomorphic and which coincide with $J^0$ outside a neighbourhood
of $V_0$. Now, if $y_0$ is the unique real ordinary cusp of $C^0$, we can proceed in the same way,
making use of the weighted homothety $\tilde{h}_\delta : (w_1 , w_2) \in B^4 \mapsto
(\delta^2 w_1 , \delta^3 w_2) \in B^4$ instead of $h_\delta$. The proof is then the same as
the one of Lemma $2.6$ of \cite{Wels1} and is not reproduced here. $\square$\\

{\bf Proof of Proposition \ref{propinteg}:}

From Lemma \ref{lemmapath}, we can assume that $J^0$ is integrable in a neighbourhood $V_0$
of the singular points of $C^0$. Let $\overline{C}_0 \subset C^0$ be a smooth compact curve
with boundary such that $C^0 \setminus \overline{C}_0 \subset V_0$. A tubular neighbourhood
$N$ of $\overline{C}_0$ in $X$ is identified with a neighbourhood of the zero section in the
normal bundle of $\overline{C}_0$ in $X$. Denote by $p : N \to \overline{C}_0$ the projection
induced by this identification and equip $N$ with the almost complex structure $J|_N$. This
identification can be chosen such that the fibres of $p$ are $J|_N$-holomorphic. 
Now as in Lemma $5.1$ of
\cite{SiTi}, there is a map $w : N \to \C$ which is holomorphic and injective once restricted
to each fibre of $p$. Such a map can be constructed as follows. Extend $J|_N$ to an almost
complex structure on the whole compactified normal bundle $\overline{N}$ over $C^0$, such that
$p : \overline{N} \to C^0$ is a sphere bundle with $J$-holomorphic fibres, the section $C^\infty$
at infinity is $J$-holomorphic and a third section $C^1$ distinct from $C^0$ and $C^\infty$ is
$J$-holomorphic. The function $w : \overline{N} \to \C P^1$ is then the unique one given by 
Riemann's uniformization theorem which is holomorphic once restricted to each fibre and
sends $C^0$, $C^\infty$ and $C^1$ to $0, \infty$ and $1$ respectively. Let $z : \overline{C}_0
\to \C$ be an injective holomorphic map, the composition with $p$ will also be denoted by
$z : N \to \C$. The antiholomorphic tangent bundle $T^{0,1}_{N , J|_N}$ is then generated by
$\partial_{\overline{w}}$ and $\partial_{\overline{z}} + a \partial_z + b \partial_w$. Moreover,
$J|_N$ is integrable if and only if $\partial_{\overline{w}} a = 0 = \partial_{\overline{w}} b$
(see \cite{SiTi}, Lemma $1.3$). This is in particular the case on $N \cap V_0$. Let 
$f : \overline{C}_0 \to \R$ be a $C^\infty$ function which is equal to one in a neighbourhood of
the boundary of $\overline{C}_0$ and to zero on $\overline{C}_0 \setminus V_0$. For 
$\delta \in [0,1]$, denote by $J_\delta$ the almost complex structure on $N$ whose
antiholomorphic tangent bundle $T^{0,1}_{N , J_\delta|_N}$ is generated by
$\partial_{\overline{w}}$ and $\partial_{\overline{z}} + \big( (1-\delta) + \delta f \circ p \big)
\big( a \partial_z + b \partial_w \big)$. Then, $J_0 = J|_N$ and $J_1$ is integrable since
$\partial_{\overline{w}} (f \circ p a) = f \circ p \partial_{\overline{w}} a = 0 $ and
$\partial_{\overline{w}} (f \circ p b) = f \circ p \partial_{\overline{w}} b = 0 $. Hence the
result. $\square$\\

Let $\gamma : t \in [0,1] \mapsto (J^t , \underline{x}^t) \in \R {\cal J}_\omega \times 
\R_{\tau} X^{c_1 (X)d - 2}$ be a path 
transversal to $\pi_\R$ and $\pi_\R^{(d_1 , m_1), (d_2 , m_2)} : 
\R {\cal M}^{(d_1 , m_1), (d_2 , m_2)}_{cusp} \to \R {\cal J}_\omega \times 
\R_{\tau_1} X^{m_1} \times \R_{\tau_2} X^{m_2}$. Let $\R {\cal M}_\gamma =
\R {\cal M}^d_{cusp} \times_{\gamma} [0,1]$, $\R \overline{\cal M}_\gamma$ its Gromov
compactification and $\pi_\gamma : \R \overline{\cal M}_\gamma \to
[0,1]$ the associated projection. 
\begin{theo}
\label{theoredcusp}
Let $\gamma : t \in [0,1] \mapsto (J^t , \underline{x}^t) \in \R {\cal J}_\omega \times 
\R_{\tau} X^{c_1 (X)d - 2}$ be a generic path chosen as above and
$C^{t_0} \in \R \overline{\cal M}_\gamma \cap \R {\cal M}^{(d_1 , m_1), (d_2 , m_2)}_{cusp}$.
Then, there exist a neighbourhood $W$ of $C^{t_0}$ in 
$\R \overline{\cal M}_\gamma$ and $\epsilon > 0$ such that $\# \pi_\gamma^{-1} (t) \cap W$ 
does not depend on the choice of $t \in ]t_0 - \epsilon , t_0 + \epsilon [ \setminus \{ t_0 \}$.
\end{theo}
Note that as soon as $W$ is small enough, all the curves of $W$ have the same mass, the one
of $C^{t_0}$.

\begin{lemma}
\label{lemmaassume}
Under the hypothesis of Theorem \ref{theoredcusp}, we can assume that $\underline{x}^t$ does
not depend on $t \in ]t_0 - \epsilon , t_0 + \epsilon [$ and that 
$(J^t)_{t \in ]t_0 - \epsilon , t_0 + \epsilon [}$ is an analytic path of almost complex structures
which are integrable in a neighbourhood of $C^{t_0}$ in $X$.
\end{lemma}

{\bf Proof:}

From Proposition \ref{propinteg}, we can assume that $J^{t_0}$ is integrable in a neighbourhood
of $C^{t_0}$. Indeed, let $J^{t_0} (\lambda)$, $\lambda \in [0,1]$, be the path in
$\R {\cal M}^{(d_1 , m_1), (d_2 , m_2)}_{cusp}$ given by this proposition, such that
$J^{t_0} (0) = J^{t_0}$ and $J^{t_0} (1)$ is integrable in a neighbourhood of $C^{t_0}$. There
is no obstruction to extend this path in a two parameters family $J^t (\lambda)$, 
$\lambda \in [0,1]$, $t \in ]t_0 - \epsilon , t_0 + \epsilon [$, such that $J^t (0) = J^t$
and $(J^t (\lambda) , \underline{x}^t)$ satisfies the hypothesis of Theorem \ref{theoredcusp}
for every $\lambda \in [0,1]$. From Lemma \ref{lemmacrit} and the definition of
$\R {\cal M}^{(d_1 , m_1), (d_2 , m_2)}_{cusp}$, none of the elements of $\pi_\gamma^{-1} (t)$
can degenerate onto a critical point of $\pi_\R$ or a reducible curve over $J^t (\lambda)$.
Thus, the cardinality of $\pi_\R^{-1} (J^t (\lambda) , \underline{x}^t) \cap W(\lambda)$,
where $W(\lambda)$ is a neighbourhood of $C^{t_0} (\lambda)$ in $\R \overline{\cal M}^d_{cusp}$, 
does not depend on $\lambda \in [0,1]$.

Now we have to prove that a transversal path $(J^t)_{t \in ]t_0 - \epsilon , t_0 + \epsilon [}$
to $\R {\cal M}^{(d_1 , m_1), (d_2 , m_2)}_{cusp}$ can be chosen analytic and made of
almost complex structures
which are integrable in a neighbourhood of $C^{t_0}$ in $X$. 
Assume that $m_1 = c_1 (X) d_1 - 1$, and
denote by $C^{t_0}_1 = [u_1, J_S^1 , J^{t_0} , \underline{x}_1 , z_c]$ the cuspidal component 
of $C^{t_0}$,
so that $\underline{x}_1 = \underline{x}^{t_0} \cap C^{t_0}_1$ is of cardinality $m_1$. 
Then, from Proposition \ref{propredrap}, the path 
$(J^t)_{t \in ]t_0 - \epsilon , t_0 + \epsilon [}$ is transversal to 
$\R {\cal M}^{(d_1 , m_1), (d_2 , m_2)}_{cusp}$ at $J^{t_0}$ if and only if
$\dot{J}^{t_0} = \frac{d}{dt} J^t|_{t = t_0}$ is such that the $J_S^1$ antilinear form
$\dot{J}^{t_0} \circ du_1 \circ J_S^1$ projects onto a generator of the cokernel
$H^1_D (S ; {\cal N}_{u_1 , -\underline{z_1}}^{cusp})_{+1} \cong \R$. Let $\alpha$ be a generator
of $H^1_D (S ; {\cal N}_{u_1 , -\underline{z_1}}^{cusp})_{+1}$ having support in a small ball $U$
of $S$. Let $V$ be an open subset of $(X , J^{t_0})$ biholomorphic to the bidisc $\overline{B}^2
\times \overline{B}^2 \subset \C^2$ and such that $V \cap C^{t_0} = u_1 (U)$. We choose this
biholomorphism such that it sends $V \cap C^{t_0}$ onto the disc $\{ w = 0 \} \subset
\overline{B}^2
\times \overline{B}^2$, where the latter is equipped with complex coordinates $(z,w)$. In this
chart, the generator $\alpha$ writes $f(z , \overline{z}) d\overline{z} \otimes w$, where
$f : \overline{B}^2 \to \C$ is $\Z/ 2\Z$-equivariant, with compact support and can be assumed
to be smooth. For every $t \in ]t_0 - \epsilon , t_0 + \epsilon [$, define then $J^t|_V$ to be
the endomorphism given by the matrix 
$\left[
\begin{array}{cc}
i & 0 \\
(t - t_0) f d\overline{z} &  i
\end{array}
\right]
$. Since $f$ is with compact support, $J^t|_V = J^{t_0}|_V$ in a neighbourhood of the boundary 
of $V$.  We can then extend $J^t$ on the whole $X$ by setting $J^t \equiv J^{t_0}$ outside $V$.
The path $t \in ]t_0 - \epsilon , t_0 + \epsilon [ \mapsto J^t$ is analytic. Moreover, the
antiholomorphic complexified tangent bundle of $(X , J^t)$ is generated by the vectors
$< \partial_{\overline{z}} + \frac{1}{2} (t-t_0) f \partial_w , \partial_{\overline{w}} >$.
Since $\frac{1}{2} (t-t_0) f$ does not depend on $w$, it follows from Lemma $1.3$ of
\cite{SiTi} that $J^t$ is integrable on $V$ for every $t \in ]t_0 - \epsilon , t_0 + \epsilon [$.
The Lemma is proved in the case $m_1 = c_1 (X) d_1 - 1$ and can be proved along the same lines
when $m_1 = c_1 (X) d_1 - 2$. $\square$\\

{\bf Proof of Theorem \ref{theoredcusp}:}

Denote by $B^2 (t_0 , \epsilon) = \{ t \in \C \, | \, |t - t_0| < \epsilon \}$. The path
$\gamma : t \in ]t_0 - \epsilon , t_0 + \epsilon [ \mapsto J^t \in \R {\cal J}_\omega$ given
by Lemma \ref{lemmaassume} is complexified to an analytic path $\gamma_\C :
t \in B^2 (t_0 , \epsilon) \mapsto J^t \in {\cal J}_\omega$ which is $\Z/2 \Z$-equivariant
and made of almost complex structures which are integrable in a neighbourhood of $C^{t_0}$.
Equip the product $Y = B^2 (t_0 , \epsilon) \times X$ with the almost complex structure
$J_Y$ defined by the matrix $\left[
\begin{array}{cc}
i & 0 \\
0 &  J^t
\end{array}
\right]
$. It is integrable in a neighbourhood of $\{ t_0 \} \times C^{t_0}$. Moreover, the sections
$(t , \underline{x}^t)$ are $J_Y$-holomorphic. Note that the complexified moduli space
${\cal M}_{\gamma_\C}$ is then a smooth curve which is equipped with a holomorphic projection
$\pi : {\cal M}_{\gamma_\C} \to B^2 (t_0 , \epsilon) \setminus \{ t_0 \}$. The complex
structure of ${\cal M}_{\gamma_\C}$ extends in a unique way on the compactification
$\overline{\cal M}_{\gamma_\C}$ and the projection 
$\pi : \overline{\cal M}_{\gamma_\C} \to B^2 (t_0 , \epsilon)$ is holomorphic. We will prove
that once restricted to any irreducible component of $\overline{\cal M}_{\gamma_\C}$, this projection
is a biholomorphism. From now on, we can assume that $\overline{\cal M}_{\gamma_\C}$ is
irreducible. Let $U \to \overline{\cal M}_{\gamma_\C}$ be the universal curve and $\overline{U}$
be the stable map compactification of $U$. The latter is a complex surface and the projection
$\overline{U} \to \overline{\cal M}_{\gamma_\C}$ is a projective line bundle having a
singular fibre over $C^{t_0} \in \overline{\cal M}_{\gamma_\C}$. Denote by $\sigma_1 , \dots ,
\sigma_{c_1(X)d - 2} : \overline{\cal M}_{\gamma_\C} \to \overline{U}$ the sections associated
to the marked points $z_1 , \dots , z_{c_1(X)d - 2}$ of $S$ and by $eval : \overline{U} \to Y$
the evaluation map, so that for every $C \in \overline{\cal M}_{\gamma_\C}$,
$eval (\sigma_j (C)) = (\pi_{\gamma_\C} (C) , x_j^{\pi_{\gamma_\C} (C)}) \in Y$. Denote by
$S_1$, $S_2$ the two irreducible components of the singular fibre of $\overline{U}$, in such a way
that $eval (S_1)$ is the cuspidal component $C^{t_0}_1$ of $C^{t_0}$. Note that since the two 
components of $C^{t_0}$
intersect transversely, $d_y eval$ is injective at the intersection point $y$ of $S_1 \cap S_2$.
The normal bundle of $S_1$ in $\overline{U}$ is isomorphic to ${\cal O}_{S_1} (-1)$ and the
evaluation map induces a non vanishing morphism from this bundle to the normal bundle
of $C^{t_0}_1$ in $Y$, which is isomorphic to ${\cal O}_{S_1} (c_1 (X) d_1 - 3) \oplus 
{\cal O}_{S_1}$. Would the projection $\pi : \overline{\cal M}_{\gamma_\C} \to B^2 (t_0 , \epsilon)$
not be biholomorphic, this morphism would vanish on $\underline{z}_1 = \underline{z} \cap S_1$
which is of cardinality $m_1$. Assume that $m_1 = c_1 (X) d_1 - 1$. The image of 
${\cal O}_{S_1} (-1)$ would then be a subline bundle of degree at least $c_1 (X) d_1 - 2$ of
${\cal O}_{S_1} (c_1 (X) d_1 - 3) \oplus 
{\cal O}_{S_1}$. This is impossible. In the same way, if $m_1 = c_1 (X) d_1 - 2$, then
the normal bundle
of $C^{t_0}_2$ in $Y$ (resp. $\overline{U}$) is isomorphic to 
${\cal O}_{S_2} (c_1 (X) d_1 - 2) \oplus {\cal O}_{S_2}$ (resp. ${\cal O}_{S_2} (-1)$).
Would the projection $\pi : \overline{\cal M}_{\gamma_\C} \to B^2 (t_0 , \epsilon)$
not be biholomorphic, the morphism ${\cal O}_{S_2} (-1) \to  {\cal O}_{S_2} (c_1 (X) d_1 - 2) 
\oplus {\cal O}_{S_2}$ would vanish on $\underline{z}_2 = \underline{z} \cap S_2$
which is of cardinality $m_2$. Since the latter is equal to $c_1 (X) d_2$, we conclude as before.
$\square$

\section{Proof of Theorem \ref{maintheo}}

Let $(J^0 , \underline{x}^0)$ and $(J^1 , \underline{x}^1)$ be two regular values of
$\pi_\R : \R {\cal M}^d_{cusp} \to \R {\cal J}_\omega \times \R_\tau X^{c_1 (X) d - 2}$ which
do not belong to $\pi_\R (\R \overline{\cal M}^d_{cusp} \setminus \R {\cal M}^d_{cusp})$.
We can assume that every real rational cuspidal $J^i$-holomorphic curve which pass through 
$\underline{x}^i$, $i \in \{ 0,1 \}$, and realize the homology class $d$ has a unique real
ordinary cusp and transversal double points as singularities, all of them being
outside $\underline{x}^i$. We have to prove that $\Gamma^d_r (J^0 , \underline{x}^0) =
\Gamma^d_r (J^1 , \underline{x}^1)$.

\subsection{Choice of the path $\gamma$}
\label{subsectchoice}

Let $\gamma : t \in [0,1] \mapsto (J^t , \underline{x}^t) \in 
\R {\cal J}_\omega \times \R_\tau X^{c_1 (X) d - 2}$ be a generic path transversal to $\pi_\R$
joining $(J^0 , \underline{x}^0)$ to $(J^1 , \underline{x}^1)$. Denote by $\R {\cal M}_\gamma =
\R {\cal M}^d_{cusp} \times_\gamma [0,1]$, $\R \overline{\cal M}_\gamma$ its Gromov
compactification and $\pi_\gamma : \R \overline{\cal M}_\gamma \to [0,1]$ the associated projection.
From Lemma \ref{lemmagencrit} and Corollary \ref{corredgen}, the path $\gamma$ can be chosen
such that all elements of $\R \overline{\cal M}_\gamma$ are irreducible real rational curves
having a unique real
ordinary cusp and only transversal double points as singularities, all of them being
outside $\underline{x}^t$, with the exception of finitely many of them which may have:

1) A unique real ordinary triple point or a unique real ordinary tacnode.

2) A real branch which crosses the real ordinary cusp transversely to its tangent line.

3) A unique real double point belonging to $\underline{x}^t$.

4) The unique real ordinary cusp belonging to $\underline{x}^t$.

5) A second real ordinary cusp outside $\underline{x}^t$. In this case, any generator
$\psi$ of $H^0 (S , K_S \otimes N_{u , -\underline{z}}^*)_{-1} 
\cong H^1 (S , N_{u , -\underline{z}})_{+1}^*$ does not vanish at cusps (compare Proposition $2.8$
of \cite{Wels1}).

6) A unique cuspidal point, which is a generic order two degenerated cusp, see \S 
\ref{subseclocalstudy} for a definition.

7) Two irreducible components.

8) Two irreducible components which intersect transversely except at one point which is a
real ordinary tacnode. These curves are not cuspidal.

In the same way, the path $\gamma$ can be chosen
such that every real rational reducible $J^t$-holomorphic curves which pass through $\underline{x}^t$
and realize the homology class $d$ have only two irreducible components, both real, and transversal double 
points as singularities, with the exception of finitely many of them which may have:

$\alpha$) A unique real ordinary triple point.

$\beta$) A unique real ordinary tacnode.

$\gamma$) A unique real ordinary cusp. In this case, any generator
$\psi$ of $H^0 (S_1 , K_S \otimes N_{u_1 , -\underline{z}_1}^*)_{-1} 
\cong H^1 (S_1 , N_{u_1 , -\underline{z}_1})_{+1}^*$ does not vanish at the cusp, where $u_1 : S_1
\to X$ parameterize the cuspidal component (compare Proposition $2.8$ of \cite{Wels1}).

$\delta$) A unique real double point belonging to $\underline{x}^t$.

$\epsilon$) Three irreducible components.

Finally, the path $\gamma$ can be chosen
such that every real rational $J^t$-holomorphic curve which pass through $\underline{x}^t$,
realize the homology class $d$, and have one of the lines $T_i$, $i \in I$, as a tangent line at
$x_i^t$, have only transversal double 
points as singularities, with the exception of finitely many which may have:

a) A unique real ordinary triple point.

b) A unique real ordinary tacnode.

c) A unique real ordinary cusp outside $\underline{x}^t$. In this case, any generator
$\psi$ of $H^0 (S , K_S \otimes N_{u , -\underline{z}}^*)_{-1} 
\cong H^1 (S , N_{u , -\underline{z}})_{+1}^*$ does not vanish at the cusp 
(compare Proposition $2.8$ of \cite{Wels1}).

d) A unique real ordinary cusp at one point $x_i^t$, $i \in I$. In this case, the curve is not 
tangent to any line $T_j$, $j \in I$. Moreover, any generator
$\psi$ of $H^0 (S , K_S \otimes N_{u , -\underline{z}}^*)_{-1} 
\cong H^1 (S , N_{u , -\underline{z}})_{+1}^*$ does not vanish at the cusp 
(compare Proposition $2.8$ of \cite{Wels1}).

e) A unique real double point belonging to $\underline{x}^t$, none of the branches being tangent
to $T_i$.

f) Two irreducible components.

g) Two irreducible components which intersect at one point $x_i^t$, $i \in I$. 
In this case, the curve is not tangent to any line $T_j$, $j \in I$.

\subsection{The case $z_c \in \protect\underline{z}$}
\label{subsectzc}

Let $(J,\underline{x}) \in \R {\cal J}_\omega \times \R_\tau X^{c_1 (X) d - 2}$ be a generic
critical value of $\pi_\R$ of type $3$ given by Lemmas \ref{lemmacrit} and \ref{lemmagencrit}.
Let $C$ be a real rational $J$-holomorphic curve which pass through 
$\underline{x}$, realize the homology class $d$ and has a unique real
ordinary cusp in $\underline{x}$, say at $x_1$. Then, the space of real rational cuspidal
$J$-holomorphic curves which pass through 
$\underline{x} \setminus \{ x_1 \}$ and realize $d$ is a one parameter family generated by
any non zero element of $H^0 (S , N_{u , -\underline{z} + z_1})_{+1} \cong \R$. In particular,
the cuspidal point moves along a smooth curve $l_c \subset \R X$ transversal to the tangent
line of $\R C$ at the cusp $x_1$. Let $U$ be a neighbourhood of  $x_1$ in $\R X$ diffeomorphic
to a ball and small enough so that $l_c$ divides it in two components. Denote by $U_-$
(resp. $U_+$) the connected component of $U \setminus l_c$ defined by the relation $U_- \cap \R C
\neq \emptyset$ (resp. $U_+ \cap \R C = \emptyset$).
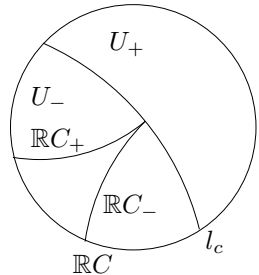
$$\vcenter{\hbox{\begin{picture}(0,0)%
\includegraphics{cusp5.pstex}%
\end{picture}%
\setlength{\unitlength}{3315sp}%
\begingroup\makeatletter\ifx\SetFigFont\undefined%
\gdef\SetFigFont#1#2#3#4#5{%
  \reset@font\fontsize{#1}{#2pt}%
  \fontfamily{#3}\fontseries{#4}\fontshape{#5}%
  \selectfont}%
\fi\endgroup%
\begin{picture}(1852,2005)(2225,-4966)
\put(3691,-4786){\makebox(0,0)[lb]{\smash{\SetFigFont{10}{12.0}{\rmdefault}{\mddefault}{\updefault}{$l_c$}%
}}}
\put(2701,-4966){\makebox(0,0)[lb]{\smash{\SetFigFont{10}{12.0}{\rmdefault}{\mddefault}{\updefault}{$\R C$}%
}}}
\put(2971,-3301){\makebox(0,0)[lb]{\smash{\SetFigFont{10}{12.0}{\rmdefault}{\mddefault}{\updefault}{$U_+$}%
}}}
\put(2926,-4516){\makebox(0,0)[lb]{\smash{\SetFigFont{10}{12.0}{\rmdefault}{\mddefault}{\updefault}{$\R C_-$}%
}}}
\put(2386,-3706){\makebox(0,0)[lb]{\smash{\SetFigFont{10}{12.0}{\rmdefault}{\mddefault}{\updefault}{$U_-$}%
}}}
\put(2386,-4021){\makebox(0,0)[lb]{\smash{\SetFigFont{10}{12.0}{\rmdefault}{\mddefault}{\updefault}{$\R C_+$}%
}}}
\end{picture}
}}$$
\begin{prop}
\label{propzc}
Let $y \in U_-$ (resp. $y \in U_+$). Then, as soon as $U$ is small enough,

1) There are exactly $2$ (resp. $0$) real rational cuspidal $J$-holomorphic curves which
pass through $\underline{x} \cup \{ y \} \setminus \{ x_1 \}$ and are close to $C$.

2) There is exactly one real rational $J$-holomorphic curves close to $C$ which
pass through $\underline{x} \cup \{ y \} \setminus \{ x_1 \}$ and has a real double point at $y$.
Moreover, this double point is real and non isolated (resp. isolated).
\end{prop}

{\bf Proof:}

As soon as $U$ is small enough, the intersection of $\R C$ with $U_-$ is made of two arcs
denoted by $\R C_+$ and $\R C_-$. The one parameter family of real rational cuspidal
$J$-holomorphic curves which pass through 
$\underline{x} \setminus \{ x_1 \}$ and realize $d$ produces one parameter families of arcs
$\R C_+$ and $\R C_-$. As soon as $U$ is small enough, these two  families provide two 
foliations of $U_-$. If $y \in U_-$, there is exactly one leaf of each foliation which pass through 
$y$, hence the first part of the proposition in this case. If $y \in U_+$, the first part of the 
proposition follows from the relation $U_+ \cap \R C = \emptyset$. Let us prove now the second
part of the proposition. First, if such a curve exists, it has to be unique. Indeed, two such
curves would intersect twice near each double point of $C$, with multiplicity one at each point of
$\underline{x} \setminus \{ x_1 \}$ and with multiplicity four at $y$. Then, their intersection
index would be at least $d^2 - c_1 (X) d + c_1 (X) d - 3 + 4 = d^2 + 1$, which is impossible.
Let $x_1' \in \R C \cap U_-$. From Proposition $2.16$ of \cite{Wels1}, the curve $\R C$ extends
to a one parameter family of real rational $J$-holomorphic curves $\R C (\eta)$, 
$\eta \in [-1 , 1]$, which pass through $\underline{x} \cup \{ x_1' \} \setminus \{ x_1 \}$
and have a non isolated (resp. isolated) real ordinary double point close to the cusp of
$\R C$ when $\eta > 0$ (resp. $\eta < 0$). In the same way as in \S \ref{subsecttac1}, the curves
$\R C (\eta)$, $\eta > 0$, form a loop in a neighbourhood of the cusp of $\R C$. Moreover, this
one parameter family of loops foliates some disk of $\R X$ centered at $x_1$. The curve $l_c$
intersects transversely these loops. Then, the non isolated real double point of
$\R C (\eta)$, $\eta > 0$, has to be in the same connected component of $U \setminus l_c$
as the branches $\R C_+$ and $\R C_-$, that is $U_-$. Thus, there exists at least one point
$y \in U_-$ for which there is a real rational $J$-holomorphic curve close to $C$ which
passes through $\underline{x} \cup \{ y \} \setminus \{ x_1 \}$ and has a non isolated
real double point at $y$. By deforming $y \in U_-$, we get the same result for all $y \in U_-$ as
soon as $U$ is small enough. The curves $\R C (\eta)$, $\eta < 0$, must then have their isolated
real double point in $U_+$, which proves the result. $\square$\\

Let us assume now that $(J, \underline{x}) = \gamma (t_0)$, where $\gamma$ is the path chosen
in \S \ref{subsectchoice} and $t_0 \in ]0,1[$. Without loss of generality, we can assume that
there exists $\epsilon > 0$ such that for every $t \in ]t_0 - \epsilon , t_0 + \epsilon[$,
$J^t \equiv J^{t_0}$ and $\underline{x}^t \setminus \{ x_1^t \} \equiv \underline{x} 
\setminus \{ x_1 \}$. The path $(x_1^t)_{t \in ]t_0 - \epsilon , t_0 + \epsilon[}$ is then
transverse to $l_c$ in $\R X$ at $t=t_0$.
\begin{prop}
\label{propzc2}
The integer $\Gamma^d_r (J^t , \underline{x}^t)$ does not depend on 
$t \in ]t_0 - \epsilon , t_0 + \epsilon[ \setminus \{ t_0 \}$.
\end{prop}

{\bf Proof:}

We can assume that for $t \in ]t_0 - \epsilon , t_0[$ (resp. $t \in ]t_0 , t_0 + \epsilon[$),
$x_1^t \in U_-$ (resp. $x_1^t \in U_+$). From the first part of Proposition \ref{propzc},
the first term in the expression of $\Gamma^d_r (J^t , \underline{x}^t)$ decreases of
$2 (-1)^{m(C^{t_0})}$ as $t$ crosses the value $t_0$. Let us prove that at the same time,
the third term in the expression of $\Gamma^d_r (J^t , \underline{x}^t)$ increases of
$2 (-1)^{m(C^{t_0})}$ as $t$ crosses the value $t_0$. From Proposition $3.4$ of
\cite{Wels1}, we have
$$\chi_{r+2}^d = \sum_{C \in {\cal T}an^d (J^t , \underline{x}^t) \, | \, T_{x_1^t} \R C = T_1} 
(-1)^{m(C)} 
+ 2 \sum_{m=0}^\delta (-1)^m \hat{n}_d^+ (m),$$
$$\chi_{r}^d = \sum_{C \in {\cal T}an^d (J^t , \underline{x}^t) \, | \, T_{x_1^t} \R C = T_1} 
(-1)^{m(C)} 
+ 2 \sum_{m=0}^\delta (-1)^m \hat{n}_d^- (m),$$
where $\chi_{r+2}^d$, $\chi_{r}^d$ are invariants and $\hat{n}_d^+ (m)$ (resp. $\hat{n}_d^- (m)$) 
is the number of real rational $J^t$-holomorphic curves of mass $m$ which pass through
$\underline{x}^t$, realize the homology class $d$ and have a non isolated (resp. isolated) real
double point at $x_1^t$. From the second part of Proposition \ref{propzc}, the term
$\sum_{m=0}^\delta (-1)^m \hat{n}_d^+ (m)$ decreases of $(-1)^{m(C^{t_0})}$ as $t$ crosses the 
value $t_0$, while the term $\sum_{m=0}^\delta (-1)^m 
\hat{n}_d^- (m)$ increases of $(-1)^{m(C^{t_0}) + 1}$. Since $\chi_{r+2}^d$, $\chi_{r}^d$ are 
constant, we deduce that 
$\sum_{C \in {\cal T}an^d (J^t , \underline{x}^t) \, | \, T_{x_1^t} \R C = T_1} (-1)^{m(C)}$
increases of
$2 (-1)^{m(C^{t_0})}$ as $t$ crosses the value $t_0$, hence the result. $\square$

\subsection{Proof of Theorem \ref{maintheo}}
\label{subsectproof}

Choose a path $\gamma : t \in [0,1] \mapsto (J^t , \underline{x}^t) \in 
\R {\cal J}_\omega \times \R_\tau X^{c_1 (X) d - 2}$ given by \S \ref{subsectchoice}.
The integer $\Gamma^d_r (J^t , \underline{x}^t)$ is then well defined for every
$t \in [0,1]$ but a finite number of parameters $0 < t_0 < \dots < t_k < 1$ corresponding
to accidents listed in \S \ref{subsectchoice}. It is obviously constant between these parameters
$t_j$, $j \in \{ 0 , \dots , k \}$, and we have to prove that is also does not change while
crossing these parameters. This is easy to check in cases $1, 2, 3, \alpha , a , b , e$ listed
in \S \ref{subsectchoice}. In cases $4 , d$, it follows from Proposition \ref{propzc2}.
The cases $5, \gamma , c$ correspond to critical points which can be treated as in
Proposition $2.16$ of \cite{Wels1}. The case $6$ follows from Proposition \ref{propcritdegcusp},
the case $7$ from Theorem \ref{theoredcusp} and cases $8, \beta$ from Theorems
\ref{theoredtac1} and \ref{theoredtac2}. Note that in this last case, the loss of two real
cuspidal curves described by Theorem \ref{theoredtac1} is compensated by the decrease of the 
multiplicity of the corresponding reducible curve. Cases $f$ and $\epsilon$ can be
treated as in Proposition $2.14$ of \cite{Wels1}. It only remains to prove the invariance
in cases $\delta , g$. While crossing such a value $t_k$, the second term in the definition
of $\Gamma^d_r (J^t , \underline{x}^t)$ remains clearly unchanged. Note that the number of
real rational $J^t$-holomorphic curves which pass through $\underline{x}^t$, 
$t \in ]t_k - \epsilon , t_k + \epsilon [ \setminus \{ t_k \}$, have $T_i$ as a tangent line
at $x_i^t$ and degenerate onto a reducible curve $C^{t_k}$ having a double point at $x_i^{t_k}$
is at most one. Indeed, two such curves would intersect with multiplicity two at $x_i^{t_k}$,
one at the other points of the configuration $\underline{x}^t$ and two near every double
point of $C^{t_k}$ but $x_i^{t_k}$. This provides a total intersection index at
least equal to $d^2 - c_1 (X)d + 2 + c_1 (X)d - 3 + 2 = d^2 + 1$, which is impossible.
Since the mass of such a real rational curve which degenerate onto such a reducible curve $C^{t_k}$
is the one of $C^{t_k}$ and the parity of the number of
real rational $J^t$-holomorphic curves which pass through $\underline{x}^t$, 
$t \in ]t_k - \epsilon , t_k + \epsilon [ \setminus \{ t_k \}$ and have $T_i$ as a tangent line
at $x_i^t$ is independant of $t$, the result follows. $\square$

\addcontentsline{toc}{part}{\hspace*{\indentation}Bibliography}


\bibliography{cuspidal}
\bibliographystyle{abbrv}

\noindent \'Ecole normale sup\'erieure de Lyon\\
Unit\'e de math\'ematiques pures et appliqu\'ees\\
UMR CNRS $5669$\\
$46$, all\'ee d'Italie\\
$69364$, Lyon cedex $07$\\
(FRANCE)\\
e-mail : {\tt jwelschi@umpa.ens-lyon.fr}

\end{document}